# A double complex for computing the sign-cohomology of the universal ordinary distribution

Greg W. Anderson
School of Mathematics
University of Minnesota
Minneapolis, MN 55455

October 21, 1997

# Contents

**Notation** The cardinality of a set $S$ is denoted $|S|$. The difference of sets $X$ and $Y$ is denoted $X \setminus Y$. The group of units of a ring $R$ is denoted $R^\times$.

# 1 Introduction

Given a positive integer $f$, a *level $f$ ordinary distribution* with values in an abelian group $A$ is a periodic function $\phi : \frac{1}{f}\mathbf{Z} \to A$ of period 1 satisfying the *level $f$ distribution relations*

$$\phi(a) = \sum_{i=0}^{g-1} \phi\left(\frac{a+i}{g}\right),$$

where $g$ is any positive integer dividing $f$ and $a \in \frac{g}{f}\mathbf{Z}$. The *universal* level $f$ ordinary distribution $U(f)$ is the quotient of the free abelian group on symbols of the form $[a]$ with $a \in \frac{1}{f}\mathbf{Z}/\mathbf{Z}$, modulo level $f$ distribution relations. An *ordinary distribution* with values in $A$ is a function $\phi : \mathbf{Q} \to A$ such that $\phi |_{\frac{1}{f}\mathbf{Z}}$ is a level $f$ ordinary distribution for all $f$; the *universal* ordinary distribution $U$ is the direct limit of the groups $U(f)$. If $\phi$ is a level $f$ ordinary distribution and $t$ is an integer prime to $f$, the function $a \mapsto \phi(ta)$ is again a level $f$ ordinary distribution, and hence $G_f := \mathrm{Gal}(\mathbf{Q}(\zeta_f)/\mathbf{Q}) = (\mathbf{Z}/f\mathbf{Z})^\times$ acts on the group $U(f)$; in the limit $G := \mathrm{Gal}(\mathbf{Q}(\zeta_\infty)/\mathbf{Q})$ acts on $U$. The terminology here is due to Kubert [10]. See Lang's book [12] for background.

Let $G_\infty \subset G$ be the subgroup generated by complex conjugation. Given any abelian group $M$ equipped with an action of $G_\infty$, we define the *sign-cohomology* (resp. *sign-homology*) of $M$ to be the Tate cohomology (resp. Tate homology) of $G_\infty$ with coefficients in $M$. The basic facts about the structure and sign-cohomology of the $G$-modules $U(f)$ and their limit $U$ are as follows. To begin with, $U(f)$ is a free abelian group of rank $|G_f|$ and the natural map $U(f) \to U$ is a split monomorphism. (In particular, the limit $U$ is a free abelian group.) Further, provided that $f > 1$ and $f \not\equiv 2 \bmod 4$, the sign-cohomology of $U(f)$ is in each degree a (free) $(\mathbf{Z}/2\mathbf{Z})$-module of rank $2^{r-1}$, where $r$ is the number of distinct primes dividing $f$, and the natural map $U(f) \to U$ induces a monomorphism in sign-cohomology. Finally, $G$ acts trivially on the sign-cohomology of $U$. These results were obtained by Kubert [10], [11].

The work of Kubert was built upon a foundation previously laid by Iwasawa and Sinnott. The $G_f$-module $U(f)$ can be regarded as a presentation by generators and relations of a $G_f$-module defined and studied much earlier by Iwasawa and subsequently employed by Sinnott [15] to compute unit- and Stickelberger-indices associated to the cyclotomic field $\mathbf{Q}(\zeta_f)$. The crucial step in Sinnott's calculation was the determination of the sign-cohomology of Iwasawa's module. In this paper we denote Iwasawa's module by $U'(f)$.

For the reader's convenience we briefly recall the definition of $U'(f)$ and sketch the construction of a $G_f$-equivariant isomorphism $U(f) \tilde{\to} U'(f)$. There exists a unique periodic function $u : \mathbf{Q} \to \mathbf{Q}$ of period 1 such that for all positive integers $f$ one has

$$\sum_{\substack{0 < x \leq f \\ (x,f)=1}} u\left(\frac{x}{f}\right) \chi(x) = \prod_{p|f} (1 - \chi(p))$$

for all primitive Dirichlet characters $\chi$ of conductor dividing $f$. Iwasawa defined $U'(f)$ to be the $\mathbf{Z}[G_f]$-submodule of $\mathbf{Q}[G_f]$ generated by elements of the form

$$\sum_{\substack{0 < x \leq f \\ (x,f)=1}} u\left(\frac{x}{g}\right) \sigma_x^{-1} \in \mathbf{Q}[G_f]$$

where $g$ is any positive integer dividing $f$ and $\sigma_x \in G_f$ satisfies $\sigma_x \zeta_f := \zeta_f^x$. Now one can verify without difficulty that $u$ is an ordinary distribution. It follows that the $G_f$-module $U'(f)$ is a quotient of $U(f)$. One can also verify easily that as a abelian group $U'(f)$ is free of rank $|G_f|$. It follows that the $G_f$-module $U'(f)$ is isomorphic to $U(f)$ because the underlying abelian groups in question are both free of rank $|G_f|$.

In this paper we study a class of problems (most of which have to do with global fields of characteristic $p > 0$) generalizing that of determining the structure and sign-cohomology of the groups $U(f)$ and their limit $U$. Our main reason for undertaking this study is to prove a conjecture of L. S. Yin [16]. In defining the generalization of $U(f)$ studied here, we more or less follow a definition given by Hayes [8] and attributed there to Mazur.

The analogue of Sinnott's unit-index calculation [15], with the Carlitz module assigned to the role played in classical cyclotomic theory by the multiplicative group, was carried out by Galovich and Rosen [6]. Quite recently, L. S. Yin [16] attempted to generalize the results of Galovich-Rosen by replacing the Carlitz module with a general sign-normalized rank one Drinfeld module. Yin computed the unit-index conditional on a remarkable conjecture concerning the Galois-module structure of the sign-cohomology of the relevant analogue of $U'(f)$. Yin's conjecture is tantalizing because it seems to be just beyond the reach of the inductive method of computation introduced by Sinnott and employed by Yin.

We prove Yin's conjecture by identifying the analogue of $U'(f)$ coming up in Yin's work with the corresponding analogue of $U(f)$, and then computing sign-cohomology of the latter by a new method involving double complexes. The machinery that we set up keeps track not only of the distribution relations themselves, but also of the higher syzygies among the distribution relations. Even in the classical cyclotomic setting our method yields a new insight: provided that $f \not\equiv 2 \bmod 4$, the sign-homology of $U(f)$ is canonically isomorphic to the Farrell-Tate homology of the subgroup of $\mathbf{Q}^\times$ generated by $-1$ and the primes dividing $f$. The Farrell-Tate homology theory, which figures prominently in our proof of Yin's conjecture, was devised by Farrell [5] to extend Tate's well known theory for finite groups to groups of finite virtual cohomological dimension. See Brown's book [1] for background. In turn, Farrell's theory has been greatly generalized by Mislin [13]. The title of the paper notwithstanding, we work with homology rather than cohomology because the former has functorial properties best suited to our purposes.

We mention that techniques developed in this paper have recently been applied by P. Das [2] to the study of *algebraic $\Gamma$-monomials*, namely complex

numbers of the form

$$\frac{\prod_i \Gamma(a_i)^{m_i}}{(2\pi i)^w}$$

where $a_i \in \mathbf{Q} \cap (0,1)$, $m_i \in \mathbf{Z}$, $w \in \mathbf{Z}$, and for all integers $t$ prime to the denominators of the $a_i$ one has

$$w = \sum_i m_i \langle ta_i \rangle$$

where $\langle x \rangle$ is the fractional part of $x$. Such numbers are in fact algebraic by a result of Koblitz and Ogus [3, Appendix] and figure in a reciprocity law due to Deligne [3],[4]; the corresponding formal sum $\sum_i m_i[a_i + \mathbf{Z}]$ represents a class in the second degree sign-cohomology of $U$ which strongly influences the Galois-theoretic properties of the monomial. Das has proved a series of results greatly illuminating the structure of the Galois group over $\mathbf{Q}$ of the extension of $\mathbf{Q}(\zeta_\infty)$ generated by the algebraic $\Gamma$-monomials. Das has also been able to give elementary proofs of some facts about algebraic $\Gamma$-monomials which previously could only be proved with the aid of Deligne's theory of absolute Hodge cycles on abelian varieties. We conclude by noting that a function field analogue of Deligne's reciprocity law recently given by S. Sinha [14] suggests that Das's theory of algebraic $\Gamma$-monomials might fruitfully be extended to global fields of characteristic $p > 0$.

## 2 Preliminaries

### 2.1 Abstract nonsense

Let $\mathcal{A}$ be an abelian category. A *chain complex* $X$ in $\mathcal{A}$ is a family $\{X_n\}_{n \in \mathbf{Z}}$ of objects of $\mathcal{A}$ equipped with a family of morphisms $\{\partial_n(X) \in \mathrm{Hom}_\mathcal{A}(X_n, X_{n-1})\}_{n \in \mathbf{Z}}$ such that $\partial_{n-1}(X)\partial_n(X) = 0$. A *chain map* $f : X \to Y$ of chain complexes in $\mathcal{A}$ is a family $\{f_n \in \mathrm{Hom}_\mathcal{A}(X_n, Y_n)\}_{n \in \mathbf{Z}}$ of morphisms such that $f_{n-1}\partial_n(X) = \partial_n(Y)f_n$. Given two chain maps $f, g : X \to Y$, a *homotopy* $T : f \to g$ is a family $\{T_n \in \mathrm{Hom}_\mathcal{A}(X_n, Y_{n+1})\}_{n \in \mathbf{Z}}$ of morphisms such that $f_n - g_n = \partial_{n+1}(Y)T_n + T_{n-1}\partial_n(X)$; we say that $f$ and $g$ are homotopic, and we write $f \sim g$, if there exists a homotopy $T : f \to g$. Given a chain complex $X$ in $\mathcal{A}$ and an integer $k$, put

$$X[k] := X_{n-k}, \quad \partial_n(X[k]) := (-1)^k \partial_{n-k}(X),$$

thereby defining the *twist* $X[k]$. Given a chain map $f : X \to Y$ of chain complexes in $\mathcal{A}$, put

$$\mathrm{Cone}(f)_n := \begin{bmatrix} X_{n-1} \\ Y_n \end{bmatrix}, \quad \partial_n(\mathrm{Cone}(f)) := \begin{bmatrix} -\partial_{n-1}(X) & 0 \\ f_{n-1} & \partial_n(Y) \end{bmatrix},$$

thereby defining the *mapping cone* $\mathrm{Cone}(f)$, which fits into a natural exact sequence

$$0 \to Y \to \mathrm{Cone}(f) \to X[1] \to 0$$

of chain complexes in $\mathcal{A}$.

**Proposition 2.1.1** *Let $f : X \to Y$ be a chain map of chain complexes in $\mathcal{A}$. Let $S$ be the set of integers $n$ such that both $f_{n-1}$ and $f_n$ are isomorphisms. Then there exists a chain map $e : \mathrm{Cone}(f) \to \mathrm{Cone}(f)$ such that $e \sim 1$ and $e_n = 0$ for all $s \in S$.*

**Proof** For each $n$, let $\phi_n : Y_n \to X_n$ be $f_n^{-1}$ or $0$ according as $f_n$ is or is not invertible. Then the family of morphisms

$$\left\{ \begin{bmatrix} 0 & \phi_n \\ 0 & 0 \end{bmatrix} : \mathrm{Cone}(f)_n \to \mathrm{Cone}(f)_{n+1} \right\}_{n \in \mathbf{Z}}$$

is a homotopy from the identity map to a map $e$ such that $e_n = 0$ for $n \in S$, as one verifies by a brief matrix calculation. ◇

**Proposition 2.1.2** *Let $g : X \to Z$ and $h : Y \to Z$ be chain maps of chain complexes in an abelian category $\mathcal{A}$. Make either of the following assumptions.*

(I) *$H_*(\mathrm{Hom}_{\mathcal{A}}(X_n, \mathrm{Cone}(h))) = 0$ for all $n$ and there exists a chain map $e : \mathrm{Cone}(h) \to \mathrm{Cone}(h)$ such that $e \sim 1$ and $e_n = 0$ for all $n \ll 0$.*

(II) *$H^*(\mathrm{Hom}_{\mathcal{A}}(X, Y_n)) = 0$ and $H^*(\mathrm{Hom}_{\mathcal{A}}(X, Z_n)) = 0$ for all $n$, and there exists a chain map $e : \mathrm{Cone}(h) \to \mathrm{Cone}(h)$ such that $e \sim 1$ and $e_n = 0$ for $n \gg 0$.*

*Then there exists a chain map $f : X \to Y$ unique up to homotopy such that $g \sim hf$.*

**Proof** Under either of the hypotheses (I) or (II), every chain map $X[k] \to \mathrm{Cone}(h)$ is homotopic to the zero map. In particular, one has $ig \sim 0$, where $i: Z \to \mathrm{Cone}(h)$ is the evident map, whence follows the existence of $f$ after a brief matrix calculation. Moreover, the difference of any two homotopies $ig \to 0$ defines a chain map $X[1] \to \mathrm{Cone}(h)$ homotopic to the zero map, whence follows the uniqueness of $f$ up to homotopy after another brief matrix calculation. ⋄

## 2.2 Farrell-Tate homology

Let $G$ be a group. We say that a (left) $G$-module $M$ (we work exclusively with left modules) is *relatively projective* if $M$ is a direct summand of a $G$-module of the form $\mathrm{Ind}_{\{1\}}^G N$ for some abelian group $N$. Here $\mathrm{Ind}_{\{1\}}^G$ is the functor left adjoint to the restriction functor $\mathrm{Res}_{\{1\}}^G$ associating to each $G$-module the underlying abelian group; more generally, given a subgroup $H \subseteq G$, the corresponding restriction functor is denoted by $\mathrm{Res}_H^G$, and the left adjoint of $\mathrm{Res}_H^G$ is denoted by $\mathrm{Ind}_H^G$.

**Proposition 2.2.1** *Let $G$ be a group, $H \subseteq G$ a subgroup of finite index, $M$ a relatively projective $G$-module, and $X$ a chain complex of $G$-modules.*
*(I) If $\mathrm{Res}_{\{1\}}^G X$ is contractible, then $H_*(\mathrm{Hom}_G(M, X)) = 0$.*
*(II) If $\mathrm{Res}_H^G X$ is contractible, then $H^*(\mathrm{Hom}_G(X, M)) = 0$.*

**Proof** There is no loss of generality in assuming that $M = \mathrm{Ind}_{\{1\}}^G N$ for some abelian group $N$.
   (I) One has $H_*(\mathrm{Hom}_G(M, X)) = H_*(\mathrm{Hom}(N, \mathrm{Res}_{\{1\}}^G X)) = 0$.
   (II) Because $H$ is of finite index in $G$, the left adjoint $\mathrm{Ind}_H^G$ and right adjoint $\mathrm{Coind}_H^G$ of the restriction functor $\mathrm{Res}_H^G$ associating to each $G$-module the underlying $H$-module are isomorphic. Therefore $H^*(\mathrm{Hom}_G(X, M)) = H^*(\mathrm{Hom}_H(\mathrm{Res}_H^G Y, \mathrm{Ind}_{\{1\}}^H N)) = 0$. ⋄

**Proposition 2.2.2** *Let $G$ be a group. Let $g: X \to Z$ and $h: Y \to Z$ be chain maps of chain complexes of $G$-modules. Assume that $X$, $Y$ and $Z$ are concentrated in nonnegative degree, $\mathrm{Res}_{\{1\}}^G \mathrm{Cone}(h)$ is contractible, and $X_n$ is relatively projective for all $n$. Then there exists a chain map $f: X \to Y$ unique up to homotopy such that $g \sim hf$.*

**Proof** This boils down to a special case of Proposition 2.1.2. ⋄

Given $G$-modules $M$ and $N$, the tensor product $M \otimes N$ is defined to be the tensor product of underlying abelian groups equipped with the diagonal $G$-action $g(m \otimes n) := (gm) \otimes (gn)$. More generally, given chain complexes $X$ and $Y$ of $G$-modules, the tensor product $X \otimes Y$ is the chain complex of $G$-modules defined by the rules

$$(X \otimes Y)_n := \bigoplus_{p+q=n} X_p \otimes Y_q$$

and

$$\partial_{p+q}(X \otimes Y)(x \otimes y) := (\partial_p(X)x) \otimes y + (-1)^p x \otimes (\partial_q(Y)y)$$

for all $x \in X_p$ and $y \in Y_q$.

We say that a chain map $f : X \to Y$ of chain complexes of $G$-modules is a *resolution* if $X$ and $Y$ are concentrated in nonnegative degree, $X_n$ is relatively projective for all $n$, and $\mathrm{Res}^G_{\{1\}}\mathrm{Cone}(f)$ is contractible. Abusing language, in a situation where the chain map $f$ is understood, we also say that $X$ is a resolution of $Y$. Proposition 2.2.2 specifies the sense in which resolutions are unique. There exists a resolution $P$ of $\mathrm{Inf}^G_{\{1\}}\mathbf{Z}$ such that $P_n$ is projective for all $n$, and consequently, given any chain complex $X$ of $G$-modules concentrated in nonnegative degree, $X \otimes P$ is a resolution of $X$. In particular, every $G$-module $M$ (viewed in this context as a chain complex of $G$-modules concentrated in degree 0) has a resolution. If $G$ is a group of cohomological dimension $r$, then there exists a resolution $P$ of $\mathrm{Inf}^G_{\{1\}}\mathbf{Z}$ such that $P_n$ is projective for all $n$ and $P_n = 0$ for $n > r$, and hence every $G$-module $M$ has a resolution $M \otimes P$ concentrated in degree $\leq r$.

**Proposition 2.2.3** *Let $G$ be a group of finite cohomological dimension. Let $P$ be a chain complex of $G$-modules such that $P_n$ is relatively projective for all $n$ and $\mathrm{Res}^G_{\{1\}}P$ is contractible. Then $P$ is contractible.*

**Proof** The complex $P^+$ obtained by replacing $P_n$ by 0 for all $n < 0$ is a resolution of $\mathrm{coker}\,\partial_0(P) = \ker \partial_{-1}(P)$, and hence has the homotopy type of a complex concentrated in degree $\leq r$, where $r$ is the cohomological dimension of $G$. It follows that $\ker \partial_n(P) = \mathrm{coker}\,\partial_{n+1}(P)$ is a direct summand of $P_n$ for all $n > r$. An evident modification of the preceding argument proves that $\ker \partial_n(P) = \mathrm{coker}\,\partial_{n+1}(P)$ is a direct summand of $P_n$ for all $n$. ⋄

**Proposition 2.2.4** *Let $G$ be a group. Let $g : X \to Z$ and $h : Y \to Z$ be chain maps of chain complexes of $G$-modules. Make the following assumptions: (i) $G$ is of finite virtual cohomological dimension. (ii) $h_n$ is an isomorphism for all $n \gg 0$. (iii) $X_n$, $Y_n$ and $Z_n$ are relatively projective for all $n$, and $\operatorname{Res}^G_{\{1\}} X$ is contractible. Then there exists a chain map $f : X \to Y$ unique up to homotopy such that $g \sim hf$.*

**Proof** By hypothesis (i) and Proposition 2.2.3 there exists a subgroup $H \subseteq G$ of finite index such that the chain complex $\operatorname{Res}^G_H X$ is contractible. By hypothesis (ii) and Proposition 2.1.1 there exists a chain map $e : \operatorname{Cone}(h) \to \operatorname{Cone}(h)$ such that $e \sim 1$ and $e_n = 0$ for all $n \gg 0$. By hypothesis (iii) and Proposition 2.2.1, one has $H^*(\operatorname{Hom}_G(X, Y_n)) = 0$ and $H^*(\operatorname{Hom}_G(X, Z_n)) = 0$ for all $n$. The result now follows by Proposition 2.1.2. $\diamond$

We say that a chain map $\kappa : X \to P$ of chain complexes of $G$-modules is a *completion* if $X_n$ and $P_n$ are relatively projective for all $n$, $\operatorname{Res}^G_{\{1\}} X$ is contractible, and $\kappa_n$ is an isomorphism for all $n \gg 0$. Proposition 2.2.4 specifies the sense in which completions are unique. Given a group $G$ of finite virtual cohomological dimension $r$, Farrell [5] (see also Brown [1, Chap. X]) showed how to construct a resolution $P$ of $\operatorname{Inf}^G_{\{1\}} \mathbf{Z}$ with $P_n$ projective for all $n$, and a completion $F \xrightarrow{\kappa} P$ with $F_n$ projective for all $n$ and $\kappa_n$ an isomorphism for all $n \geq r$; the tensor product $M \otimes P$ is then a resolution of $M$, and the tensor product $M \otimes F$ a completion of $M \otimes P$.

Given a group $G$ of finite virtual cohomological dimension and a $G$-module $M$, one defines
$$\hat{H}_*(G, M) := H_*(\operatorname{Coinv}^G_G X)$$
where $X$ is any completion of a resolution of $M$, and $\operatorname{Coinv}^G_G$ is the functor left adjoint to the functor $\operatorname{Inf}^G_{\{1\}}$ equipping abelian groups with trivial $G$-action. We also introduce the abbreviated notation
$$\hat{H}_*(G) := \hat{H}_*(G, \operatorname{Inf}^G_{\{1\}} \mathbf{Z}).$$

The *Farrell-Tate homology theory* $\hat{H}_*$ extends to groups of finite virtual cohomological dimension the theory introduced by Tate for finite groups.

## 2.3 The Shapiro lemma and related results

Let $G$ be a group of finite virtual cohomological dimension, and let $H$ be a subgroup (necessarily also of finite virtual cohomological dimension). Let $N$ be an $H$-module and let $Y$ be a completion of a resolution of $N$. Then $\operatorname{Ind}_H^G Y$ is a completion of a resolution of $\operatorname{Ind}_H^G N$. Further, the functors $\operatorname{Res}_H^G \circ \operatorname{Inf}_{\{1\}}^G$ and $\operatorname{Inf}_{\{1\}}^H$ are isomorphic, and hence so are their left adjoints $\operatorname{Coinv}_G^G \circ \operatorname{Ind}_H^G$ and $\operatorname{Coinv}_H^H$. One thus obtains canonical isomorphisms

$$\begin{aligned}
\hat{H}_*(G, \operatorname{Ind}_H^G N) &= H_*(\operatorname{Coinv}_G^G \operatorname{Ind}_H^G X)) \\
&= H_*(\operatorname{Coinv}_H^H X) \\
&= \hat{H}_*(H, N)
\end{aligned}$$

of graded abelian groups. The assertion that there exists an isomorphism between the extreme terms in the relation above, functorial in $H$-modules $N$, is the *Shapiro lemma* for Farrell-Tate homology.

With $G$ and $H$ as in the preceding paragraph, let $M$ be a $G$-module and let $X$ be a completion of a resolution of $M$. Then $\operatorname{Res}_H^G X$ is a completion of a resolution of $\operatorname{Res}_H^G M$. Suppose now that $H$ is a normal subgroup of $G$ and put $Q := G/H$. Let $\operatorname{Coinv}_H^G$ be the functor left adjoint to the inflation functor $\operatorname{Inf}_Q^G$. The functors $\operatorname{Coinv}_H^H \circ \operatorname{Res}_H^G$ and $\operatorname{Res}_{\{1\}}^Q \circ \operatorname{Coinv}_H^G$ are isomorphic. Moreover, the functor $\operatorname{Res}_{\{1\}}^Q$ is exact. One thus obtains canonical isomorphisms

$$\begin{aligned}
\hat{H}_*(H, \operatorname{Res}_H^G M) &= H_*(\operatorname{Coinv}_H^H \operatorname{Res}_H^G X) \\
&= H_*(\operatorname{Res}_{\{1\}}^Q \operatorname{Coinv}_H^G X) \\
&= \operatorname{Res}_{\{1\}}^Q H_*(\operatorname{Coinv}_H^G X)
\end{aligned}$$

of graded abelian groups. Thus $\hat{H}_*(H, \operatorname{Res}_H^G M)$ is canonically equipped with graded $Q$-module structure; in the sequel we identify $\hat{H}_*(H, \operatorname{Res}_H^G M)$ with $H_*(\operatorname{Coinv}_H^G X)$ rather than $H_*(\operatorname{Coinv}_H^H \operatorname{Res}_H^G X)$.

**Proposition 2.3.1** *Let $G$ be a group of finite virtual cohomological dimension and let $H \subseteq G$ be a normal subgroup. Let $\sigma$ be an element of the center of $G$. Let $M$ be a $G$-module on which $\sigma$ acts trivially. Let $X$ be a completion of a resolution of $M$. Then the automorphism of $X$ induced by $\sigma$ is homotopic to the identity, and hence $\sigma H \in G/H$ induces the identity mapping in $\hat{H}_*(H, Res_H^G M) = H_*(Coinv_H^G X)$.*

**Proof** This is a consequence of the uniqueness of completions of resolutions (Proposition 2.2.2 and Proposition 2.2.4). ⋄

**Proposition 2.3.2** *Let $\Gamma$ be a group of finite virtual cohomological dimension. Let $G$ be a normal subgroup of $\Gamma$. Let $\Pi$ be any subgroup of $\Gamma$. Put*

$$H := G \cap \Pi, \quad \bar{\Gamma} := \Gamma/G, \quad \bar{\Pi} := \Pi/H.$$

*Let $M$ be a $\Pi$-module. Then there exists an isomorphism*

$$\hat{H}_*(G, Res_G^\Gamma Ind_\Pi^\Gamma M) = Ind_{\bar{\Pi}}^{\bar{\Gamma}} \hat{H}_*(H, Res_H^\Pi M)$$

*of $\bar{\Gamma}$-modules functorial in $M$.*

**Proof** Clearly the functors $\mathrm{Res}_\Pi^\Gamma \circ \mathrm{Inf}_{\bar{\Gamma}}^\Gamma$ and $\mathrm{Inf}_{\bar{\Pi}}^\Pi \circ \mathrm{Res}_{\bar{\Pi}}^{\bar{\Gamma}}$ are isomorphic, and hence so are their left adjoints $\mathrm{Coinv}_G^\Gamma \circ \mathrm{Ind}_\Pi^\Gamma$ and $\mathrm{Ind}_{\bar{\Pi}}^{\bar{\Gamma}} \circ \mathrm{Coinv}_H^\Pi$. Moreover, the functor $\mathrm{Ind}_{\bar{\Pi}}^{\bar{\Gamma}}$ is exact. Let $X$ be a completion of a resolution of $M$. One has canonical isomorphisms

$$\begin{aligned}
\hat{H}_*(G, \mathrm{Res}_G^\Gamma \mathrm{Ind}_\Pi^\Gamma M) &= H_*(\mathrm{Coinv}_G^\Gamma \mathrm{Ind}_\Pi^\Gamma X) \\
&= H_*(\mathrm{Ind}_{\bar{\Pi}}^{\bar{\Gamma}} \mathrm{Coinv}_H^\Pi X) \\
&= \mathrm{Ind}_{\bar{\Pi}}^{\bar{\Gamma}} H_*(\mathrm{Coinv}_H^\Pi X) \\
&= \mathrm{Ind}_{\bar{\Pi}}^{\bar{\Gamma}} \hat{H}_*(H, \mathrm{Res}_H^\Pi M)
\end{aligned}$$

of graded $\bar{\Gamma}$-modules. ⋄

## 2.4 The double complex $\mathcal{KT}$

We give a construction exploited repeatedly in the paper. The input for the construction is as follows:

- A commutative ring $R$ with unit.
- An $R$-module $M$.
- A linearly ordered set $S$.

- A family $\{f_s \in R\}_{s \in S}$.

- Elements $f^\pm$ of $R$ such that $f^+ f^- = 0$.

The output of the construction is as follows:

- A double complex $\mathcal{KT}$ of $R$-modules, i. e., a bigraded $R$-module

$$\mathcal{KT} = \mathcal{KT}\left(M/R, \{f_s\}_{s \in S}, \begin{bmatrix} f^+ \\ f^- \end{bmatrix}\right) = \bigoplus_m \bigoplus_n \mathcal{KT}_{mn}$$

equipped with $R$-linear maps

$$\partial, \delta : \mathcal{KT} \to \mathcal{KT}$$

of bidegree $(-1, 0)$ and $(0, -1)$, respectively, such that $\partial^2 = 0$, $\delta^2 = 0$, and $\partial \delta + \delta \partial = 0$.

- Chain complexes of $R$-modules

$$\begin{aligned}
\mathcal{K} &= \mathcal{K}\left(M/R, \{f_s\}_{s \in S}\right), \\
\overline{\mathcal{K}} &= \overline{\mathcal{K}}\left(M/R, \{f_s\}_{s \in S}, \begin{bmatrix} f^+ \\ f^- \end{bmatrix}\right), \\
\mathcal{T} &= \mathcal{T}\left(M/R, \begin{bmatrix} f^+ \\ f^- \end{bmatrix}\right), \\
\overline{\mathcal{T}} &= \overline{\mathcal{T}}\left(M/R, \{f_s\}_{s \in S}, \begin{bmatrix} f^+ \\ f^- \end{bmatrix}\right), \\
\mathcal{KT}^{tot} &= \mathcal{KT}^{tot}\left(M/R, \{f_s\}_{s \in S}, \begin{bmatrix} f^+ \\ f^- \end{bmatrix}\right), \\
\mathcal{KT}^+ &= \mathcal{KT}^+\left(M/R, \{f_s\}_{s \in S}, \begin{bmatrix} f^+ \\ f^- \end{bmatrix}\right),
\end{aligned}$$

which we call the *companions* of the double complex $\mathcal{KT}$.

The notation $\mathcal{KT}$ is meant to call Koszul and Tate to mind.

Here is the construction. We define $\mathcal{S}$ to be the free abelian group on symbols of the form $[I, k]$ where $I \subseteq S$ is a finite subset and $k$ is an integer, and we bigrade $\mathcal{S}$ by declaring the symbol $[I, k]$ to be of bidegree $(|I|, k)$. Put

$$\mathcal{KT} := M \otimes \mathcal{S}$$

$$\partial(m \otimes [I, k]) := \sum_{i \in I} (-1)^{|\{j \in I | j < i\}|} f_i m \otimes [I \setminus \{i\}, k]$$

$$\delta(m \otimes [I, k]) := (-1)^{|I|} \begin{cases} f^+ m \otimes [I, k-1] & \text{if } k \text{ is even} \\ f^- m \otimes [I, k-1] & \text{if } k \text{ is odd} \end{cases}$$

for all $m \in M$, finite subsets $I \subseteq S$ and integers $k$. Take $\mathcal{K}_n := \mathcal{KT}_{n0}$ and equip $\mathcal{K}$ with the differential induced by the operator $\partial$. Take $\mathcal{T}_n := \mathcal{KT}_{0n}$ and equip $\mathcal{T}$ with the differential induced by the operator $\delta$. Put

$$\overline{\mathcal{K}} := \mathcal{K}\left(\frac{M}{f^- M}/R, \{f_s\}_{s \in S}\right), \quad \overline{\mathcal{T}} := \mathcal{T}\left(\frac{M}{\sum_{s \in S} f_s M}/R, \begin{bmatrix} f^+ \\ f^- \end{bmatrix}\right).$$

Let $\mathcal{KT}^{tot}$ be the *total complex* associated to the double complex $\mathcal{KT}$, i. e., a copy of $\mathcal{KT}$ graded by total degree and equipped with the differential induced by $\partial + \delta$. Let $\mathcal{KT}^-$ be the subcomplex of $\mathcal{KT}^{tot}$ spanned over $R$ by elements of the form $m \otimes [I, k]$ with $m \in M$, $I \subseteq S$ finite, and $k < 0$. Finally, put $\mathcal{KT}^+ := \mathcal{KT}^{tot}/\mathcal{KT}^-$. Note that $\overline{\mathcal{T}}$ is naturally a quotient of $\mathcal{KT}^{tot}$ and $\overline{\mathcal{K}}$ naturally a quotient of $\mathcal{KT}^+$.

**Proposition 2.4.1** *Let $R$ be a commutative ring with unit, $M$ an $R$-module, $\{f_s\}_{s \in S}$ a family of elements of $R$ indexed by a linearly ordered set $S$, $f^\pm$ elements of $R$ such that $f^+ f^- = 0$. Consider the double complex*

$$\mathcal{KT}\left(M/R, \{f_s\}_{s \in S}, \begin{bmatrix} f^+ \\ f^- \end{bmatrix}\right)$$

*and companion complexes $\mathcal{K}$, $\overline{\mathcal{K}}$, $\mathcal{T}$, $\overline{\mathcal{T}}$, $\mathcal{KT}^{tot}$ and $\mathcal{KT}^+$.*

*(I) If, for all finite subsets $I \subseteq S$, the sequence $\{f_i\}_{i \in I}$ is $M$-regular, then $\mathcal{K}$ is acyclic in positive degree.*

*(II) If, for some $s \in S$, the element $f_s$ operates invertibly on $M$, then $\mathcal{K}$ is acyclic.*

(III) *If $\mathcal{K}$ is acyclic in positive degree, then the quotient map $\mathcal{KT}^{tot} \to \overline{\mathcal{T}}$ induces an isomorphism in homology.*

(IV) *If $\mathcal{T}$ is acyclic, then $\mathcal{KT}^{tot}$ is acyclic and the quotient map $\mathcal{KT}^+ \to \overline{\mathcal{K}}$ induces an isomorphism in homology.*

**Proof** Because homology commutes with direct limits, we may assume that $S$ is a finite set. Then assertions (I) and (II) are standard facts about Koszul complexes; assertions (III) and (IV) are proved by straightforward spectral sequence arguments. ⋄

## 2.5 Almost free abelian groups

Finitely generated abelian groups are of finite virtual cohomological dimension and hence the Farrell-Tate theory applies to them.

**Proposition 2.5.1** *Let $\Gamma$ be a finitely generated abelian group and let $G \subset \Gamma$ be a subgroup. Let $\Delta \subseteq G \times G$ be the diagonal subgroup. Let $p : \Gamma \times G \to \Gamma$ and $q : \Gamma \times G \to G$ be the first and second projections, respectively. Let $r : \Gamma \xrightarrow{\sim} (\Gamma \times G)/\Delta$ be the isomorphism inverse to that induced by $p - q$. Let $F$ be a completion of a resolution $P$ of $\mathrm{Inf}_{\{1\}}^G \mathbf{Z}$ such that $F_n$ and $P_n$ projective for all $n$. Let $M$ be a $\Gamma$-module. Put*

$$M' := r^* \mathrm{Coinv}_{\Delta}^{\Gamma \times G}(p^*M \otimes q^*F).$$

*Then there exists an isomorphism*

$$H_*(M') = \mathrm{Inf}_{\Gamma/G}^{\Gamma} \hat{H}_*(G, \mathrm{Res}_G^{\Gamma} M)$$

*of graded $\Gamma$-modules functorial in $M$.*

**Proof** Without loss of generality we may assume that $F = \mathrm{Res}_G^{\Gamma} \tilde{F}$, where $\tilde{F}$ is a completion of a resolution $\tilde{P}$ of $\mathrm{Inf}_{\{1\}}^{\Gamma} \mathbf{Z}$ such that $\tilde{P}_n$ and $\tilde{F}_n$ are projective for all $n$. Let

$$\left. \begin{array}{l} i := (\gamma \mapsto (\gamma, 1)) \\ d := (\gamma \mapsto (\gamma, \gamma)) \end{array} \right\} : \Gamma \to \Gamma \times \Gamma$$

and let $\tilde{q} : \Gamma \times \Gamma \to \Gamma$ be the second projection. Consider the complex

$$\tilde{M} := \mathrm{Inf}_{(\Gamma \times \Gamma)/\Delta}^{\Gamma \times \Gamma} \mathrm{Coinv}_{\Delta}^{\Gamma \times \Gamma}(p^*M \otimes \tilde{q}^* \tilde{F})$$

of $\Gamma \times \Gamma$-modules. Now $d^*\tilde{M} = M \otimes \tilde{F}$ is a completion of a resolution of $M$ and one has an isomorphism

$$d^* \circ \operatorname{Inf}^{\Gamma \times \Gamma}_{(\Gamma \times \Gamma)/\Delta} \circ \operatorname{Coinv}^{\Gamma \times \Gamma}_{\Delta} = \operatorname{Inf}^{\Gamma}_{\Gamma/G} \circ \operatorname{Coinv}^{\Gamma}_{G} \circ d^*$$

of functors; accordingly, we have an isomorphism

$$H_*(d^*\tilde{M}) = \operatorname{Inf}^{\Gamma}_{\Gamma/G} \hat{H}_*(G, \operatorname{Res}^{\Gamma}_{G} M)$$

of graded $\Gamma$-modules functorial in $M$. One has an isomorphism of functors

$$i^* \operatorname{Inf}^{\Gamma \times \Gamma}_{(\Gamma \times \Gamma)/\Delta} \operatorname{Coinv}^{\Gamma \times \Gamma}_{\Delta} = r^* \operatorname{Coinv}^{\Gamma \times G}_{\Delta} \operatorname{Res}^{\Gamma \times \Gamma}_{\Gamma \times G}$$

and thus we have an isomorphism

$$H_*(i^*\tilde{M}) = H_*(M')$$

of graded $\Gamma$-modules functorial in $M$. Finally, for all $\gamma \in \Gamma$, the elements $i(\gamma)$ and $d(\gamma)$ of $\Gamma \times \Gamma$ induce homotopic automorphisms of the complex $\tilde{M}$ by Proposition 2.3.1, and hence we have an isomorphism

$$H_*(d^*\tilde{M}) = H_*(i^*\tilde{M})$$

of graded $\Gamma$-modules functorial in $M$. ◇

We say that an abelian group is *almost free* if the group can be factored as the product of a free abelian group and a finite cyclic group. The multiplicative group of a global field is almost free. Every subgroup of an almost free abelian group is again almost free.

**Proposition 2.5.2** *Let $\Gamma$ be a finitely generated abelian group. Let $G \subseteq \Gamma$ be an almost free subgroup of rank $r$, and let $m$ be the order of the torsion subgroup of $G$. Let $g_1, \ldots, g_r \in G$ be independent, and let $g_0 \in G$ generate the torsion subgroup of $G$. Let $M$ be a $\Gamma$-module. Consider the chain complex*

$$M' := \mathcal{KT}^{tot}\left(M/\mathbf{Z}[\Gamma], \{1 - g_i\}_{i=1}^{r}, \begin{bmatrix} \sum_{i=0}^{m-1} g_0^i \\ 1 - g_0 \end{bmatrix}\right).$$

*of $\Gamma$-modules. Then there exists an isomorphism*

$$H_*(M') = \operatorname{Inf}^{\Gamma}_{\Gamma/G} \hat{H}_*(G, \operatorname{Res}^{\Gamma}_{G} M)$$

*of graded $\Gamma$-modules functorial in $M$.*

**Proof** Consider the chain complexes

$$P := \mathcal{KT}^{+}\left(\mathbf{Z}[G]/\mathbf{Z}[G], \{1-g_i^{-1}\}_{i=1}^r, \begin{bmatrix} \sum_{i=0}^{m-1} g_0^{-i} \\ 1-g_0^{-1} \end{bmatrix}\right)$$

and

$$F := \mathcal{KT}^{tot}\left(\mathbf{Z}[G]/\mathbf{Z}[G], \{1-g_i^{-1}\}_{i=1}^r, \begin{bmatrix} \sum_{i=0}^{m-1} g_0^{-i} \\ 1-g_0^{-1} \end{bmatrix}\right)$$

of $G$-modules. Then $P_n$ and $F_n$ are projective for all $n$ and moreover $F_n = P_n$ for all $n \geq r$; by Proposition 2.4.1, $P_n$ is a resolution of $\mathrm{Inf}_{\{1\}}^G \mathbf{Z}$ and $F$ a completion of $P$. The result is now just a special case of Proposition 2.5.1. ⋄

**Proposition 2.5.3** *Let $G$ be an almost free abelian group of positive finite rank $r$ and let $m$ be the order of the torsion subgroup of $G$. Then $H_n(G)$ is a free $(\mathbf{Z}/m\mathbf{Z})$-module of rank $2^{r-1}$ for all $n$.*

**Proof** By Proposition 2.5.2 we have an isomorphism

$$H_n(G)\tilde{\to} H_n\left(\mathcal{KT}^{tot}\left(\mathbf{Z}/\mathbf{Z}, \{0\}_{i=1}^r, \begin{bmatrix} m \\ 0 \end{bmatrix}\right)\right),$$

whence the result after a brief computation with binomial coefficients. ⋄

**Proposition 2.5.4** *Let $\Gamma$ be an abelian group (not necessarily finitely generated) and let $\Pi \subseteq \Gamma$ be a subgroup of finite index. Let $G \subseteq \Gamma$ be an almost free subgroup of finite rank $r$ and let $m$ be the order of the torsion subgroup of $G$. Let $g_1, \ldots, g_r \in G$ be independent and let $g_0 \in G$ generate the torsion subgroup. Consider the chain complex*

$$K := \mathcal{KT}^{tot}\left(\mathrm{Ind}_\Pi^\Gamma \mathrm{Inf}_{\{1\}}^\Pi \mathbf{Z}/\mathbf{Z}[\Gamma], \{1-g_i\}_{i=1}^r, \begin{bmatrix} \sum_{i=0}^{m-1} g_0^i \\ 1-g_0 \end{bmatrix}\right)$$

*of $\Gamma$-modules. Then one has an isomorphism*

$$H_*(K) \tilde{\to} \mathrm{Ind}_{\Pi G}^\Gamma \mathrm{Inf}_{\{1\}}^{\Pi G} \hat{H}_*(\Pi \cap G)$$

*of graded $\Gamma$-modules.*

**Proof** Let $\Gamma' \subseteq \Gamma$ be a finitely generated subgroup such that $G \subseteq \Gamma'$ and $\Gamma'\Pi = \Gamma$. Replacing $\Gamma$ with $\Gamma'$, we may assume that $\Gamma$ is finitely generated. We have at our disposal a canonical isomorphism

$$H_*(K) = \mathrm{Inf}_{\Gamma/G}^{\Gamma} \hat{H}_*(G, \mathrm{Res}_G^{\Gamma} \mathrm{Ind}_{\Pi}^{\Gamma} \mathrm{Inf}_{\{1\}}^{\Pi} \mathbf{Z})$$

provided by Proposition 2.5.2, a canonical isomorphism

$$\hat{H}_*(G, \mathrm{Ind}_{\Pi}^{\Gamma} \mathrm{Inf}_{\{1\}}^{\Pi} \mathbf{Z}) = \mathrm{Ind}_{\Pi/(\Pi \cap G)}^{\Gamma/G} \hat{H}_*(\Pi \cap G, \mathrm{Res}_{\Pi \cap G}^{\Pi} \mathrm{Inf}_{\{1\}}^{\Pi} \mathbf{Z})$$

provided by Proposition 2.3.2, a canonical isomorphism

$$\hat{H}_*(\Pi \cap G, \mathrm{Res}_{\Pi \cap G}^{\Pi} \mathrm{Inf}_{\{1\}}^{\Pi} \mathbf{Z}) = \mathrm{Inf}_{\{1\}}^{\Pi \cap (\Pi \cap G)} \hat{H}_*(\Pi \cap G)$$

provided by Proposition 2.3.1, and an isomorphism

$$\mathrm{Inf}_{\Gamma/G}^{\Gamma} \circ \mathrm{Ind}_{\Pi/(\Pi \cap G)}^{\Gamma/G} \circ \mathrm{Inf}_{\{1\}}^{\Pi/(\Pi \cap G)} = \mathrm{Ind}_{\Pi G}^{\Gamma} \circ \mathrm{Inf}_{\{1\}}^{\Pi G}$$

of functors, whence the result. ◇

# 3 The principal objects of study

## 3.1 The basic data $(\mathbf{K}, A, \mathrm{sgn})$

For the rest of the paper we fix the following items.

- A locally compact nondiscrete topological field $\mathbf{K}$ containing only finitely many roots of unity.

- A discrete cocompact integrally closed subring $A \subset \mathbf{K}$.

- A continuous homomorphism

  $$\mathrm{sgn} : \mathbf{K}^{\times} \to \text{(the group of roots of unity in } \mathbf{K}\text{)}$$

  the restriction of which to the group of roots of unity of $\mathbf{K}$ is the identity mapping.

We call sgn the *sign homomorphism*, and we say that $x \in \mathbf{K}^\times$ is *positive* if sgn $x = 1$. We say that the basic data are *archimedean* if $\mathbf{K}$ is archimedean. We denote the fraction field of $A$ by $k$.

There is only one archimedean example $(\mathbf{K}, A, \text{sgn})$ of basic data, namely the triple $(\mathbf{R}, \mathbf{Z}, x \mapsto x/|x|)$. In the archimedean case our usage of the term "positive" is just the ordinary usage.

The simplest example $(\mathbf{K}, A, \text{sgn})$ of nonarchimedean basic data arises as follows. Let $\mathbf{F}_q$ be the field of $q$ elements and let $\mathbf{F}_q(T)$ be the field of rational functions in a variable $T$ with coefficients in $\mathbf{F}_q$. Take $\mathbf{K}$ to be the completion $\mathbf{F}_q((1/T))$ of $\mathbf{F}_q(T)$ at the infinite place. Take $A$ to be the polynomial ring $\mathbf{F}_q[T]$. Decompose $\mathbf{F}_q((1/T))^\times$ as a direct product

$$\mathbf{F}_q^\times \cdot T^{\mathbf{Z}} \cdot (1 + (1/T)\mathbf{F}_q[[1/T]])$$

and take the sign homomorphism sgn to be projection to the first factor. In this example the positive elements of $\mathbf{F}_q[T]$ are the monic polynomials.

Every nonarchimedean example $(\mathbf{K}, A, \text{sgn})$ of basic data arises in the following manner. Let $X/k_0$ be a smooth projective geometrically irreducible curve defined over a finite field $k_0$. Let $\infty$ be a closed point of $X$. Take $\mathbf{K}$ to be the completion of the function field of $X$ at $\infty$. Take $A$ to be the ring consisting of elements of the function field of $X$ regular away from $\infty$. Choose an isomorphism of topological fields under which to identify $\mathbf{K}$ with $\mathbf{F}_q((1/T))$, where $q$ is the cardinality of the residue field of $\infty$, and define the sign homomorphism as in the preceding example.

For archimedean and nonarchimedean basic data $(\mathbf{K}, A, \text{sgn})$ alike, the ring $A$ is a Dedekind domain the group of units of which is finite, the class group of which is finite, and every residue field of which is finite.

## 3.2 $A$-lattices

A *fractional $A$-ideal* is a finitely generated nonzero $A$-submodule of $k$. An *integral $A$-ideal* is a fractional $A$-ideal contained in $A$. When we speak of fractional or integral $A$-ideals we tacitly exclude the zero ideal of $A$ from consideration. We say that a fractional $A$-ideal $I$ is *principal in the narrow sense* if $I = (a)$ for some positive $a \in k^\times$. The quotient of the group of fractional $A$-ideals by the subgroup of ideals principal in the narrow sense is by definition the *narrow ideal class group*. The narrow ideal class group is finite.

An *A-lattice* is by definition a cocompact discrete $A$-submodule of $\mathbf{K}$. An $A$-lattice is without $A$-torsion and contains a copy of $A$ as a subgroup of finite index, and therefore is projective over $A$ of rank one. We say that two $A$-lattices $W_1$ and $W_2$ are *homothetic*, and we write $W_1 \sim W_2$, if there exists some positive $x \in \mathbf{K}$ such that $xW_1 = W_2$. The relation of homothety is an equivalence relation in the set of $A$-lattices.

Every fractional $A$-ideal is an $A$-lattice. Fractional $A$-ideals belong to the same narrow ideal class if and only if they are homothetic. Every homothety class of $A$-lattices contains a fractional $A$-ideal. The set of homothety classes of $A$-lattices thus corresponds bijectively with the narrow ideal class group and in particular is finite.

## 3.3 The set $\Xi$

Given $x \in \mathbf{K}$ and an $A$-lattice $W$, we say that

$$x + W := \{x + w \in \mathbf{K} \mid w \in W\}$$

is a *translate* of the $A$-lattice $W$, and we say that a subset of $\mathbf{K}$ is an *A-lattice translate* if of the form $x + W$. We always write $A$-lattice translates as a sum, the first symbol denoting an element of $\mathbf{K}$ and the second an $A$-lattice. We say that two $A$-lattice translates $x_1 + W_1$ and $x_2 + W_2$ are *homothetic*, and we write $x_1 + W_1 \sim x_2 + W_2$, if for some positive $y \in \mathbf{K}$ one has $yW_1 = W_2$ and $yx_1 - x_2 \in W_2$. Homothety is an equivalence relation in the set of $A$-lattice translates. Note that for all $A$-lattices $W$ and $x, y \in \mathbf{K}$, one has $x + W \sim y + W$ if and only if $x - y \in W$. Given an $A$-lattice translate $x + W$ and an integral $A$-ideal $f$, we say that $x + W$ is *annihilated by $f$*, or that $x + W$ is *$f$-torsion*, if $xf \subseteq W$; and we say that $x + W$ is of *order $f$* if $\{a \in A \mid ax \in W\} = f$. We say that an $A$-lattice translate is *torsion* if $f$-torsion for some integral $A$-ideal $f$.

We denote the homothety class of a torsion $A$-lattice translate $x + W$ by $[x+W]$. We denote the set of homothety classes of torsion $A$-lattice translates by $\Xi$. For each integral $A$-ideal $f$, let $\Xi(f)$ be the set of homothety classes of $f$-torsion $A$-lattice translates, and let $\Xi^\times(f)$ be the set of homothety classes of $A$-lattice translates of order $f$. Given a fractional $A$-ideal $I$ and an $A$-lattice $W$ let $I \cdot W$ be the $A$-submodule of $\mathbf{K}$ generated by all products of the form $aw$ where $a \in I$ and $w \in W$; the $A$-submodule $I \cdot W$ is again an $A$-lattice.

**Proposition 3.3.1** *Let $f$ be an integral $A$-ideal. There exists a unique map*

$$Y_f : \Xi \to \Xi$$

*such that*

$$Y_f[x + W] = [x + f^{-1} \cdot W]$$

*for all torsion $A$-lattice translates $x + W$. Every fiber of the map $Y_f$ is of cardinality $|A/f|$. One has*

$$Y_f \circ Y_g = Y_{fg}$$

$$Y_f \Xi(g) = \begin{cases} \Xi(g/f) & \text{if } f \text{ divides } g \\ \Xi(g) & \text{if } f \text{ is prime to } g \end{cases}$$

$$Y_f \Xi^\times(g) = \begin{cases} \Xi^\times(g/f) & \text{if } f \text{ divides } g \\ \Xi^\times(g) & \text{if } f \text{ is prime to } g \end{cases}$$

$$Y_f^{-1} \Xi(g) = \Xi(fg)$$

*for all integral $A$-ideals $g$.*

**Proof** The proof is quite straightforward and we omit it. ⋄

## 3.4 The profinite group $G$

Given integral $A$-ideals $f$, $I$ and $J$, we write $I \sim_f J$ if $I$ and $J$ are prime to $f$ and there exist nonzero $a, b \in A$ prime to $f$ such that $b/a$ is positive, $a \equiv b \bmod f$ and $aI = bJ$. The relation $\sim_f$ is an equivalence relation in the set of integral $A$-ideals prime to $f$. For each integral $A$-ideal $f$, the quotient $G_f$ of the monoid of integral $A$-ideals prime to $f$ by the equivalence relation $\sim_f$ is a finite abelian group. The family of groups $\{G_f\}$ forms an inverse system indexed by the set of integral $A$-ideals directed by the divisibility relation. Put

$$G := \varprojlim G_f.$$

The transition maps of the inverse system $\{G_f\}$ are surjective and hence each group $G_f$ is canonically a quotient of $G$.

Given an integral $A$-ideal $f$, let $\sqrt{f}$ be the product of the maximal $A$-ideals dividing $f$.

**Proposition 3.4.1** *For each integral A-ideal $f$, the natural map*

$$\prod_{p|f} \ker\left(G_f \to G_{f/p}\right) \to \ker\left(G_f \to G_{f/\sqrt{f}}\right)$$

*is bijective, where the Cartesian product on the left is extended over the maximal A-ideals $p$ dividing $f$.*

**Proof** There exists a unique isomorphism

$$(A/f)^\times \xrightarrow{\sim} \ker(G_f \to G_1)$$

under which, for all positive $a \in A$ prime to $f$, the congruence class of $a$ modulo $f$ maps to the $\sim_f$-equivalence class of $(a)$. This noted, the proposition reduces to an instance of the Chinese Remainder Theorem. ◇

**Proposition 3.4.2** *Let $f$ be an integral A-ideal. Let $x + W$ be a torsion A-lattice translate of order $f$.*

*(I) For all integral A-ideals $I$ and $J$ both prime to $f$, one has $I \sim_f J$ if and only if $x + I^{-1} \cdot W \sim x + J^{-1} \cdot W$.*

*(II) For all torsion A-lattice translates $x' + W'$ of order $f$ there exists an integral A-ideal $I$ prime to $f$ such that $x' + W' \sim x + I^{-1} \cdot W$.*

**Proof** (I)($\Rightarrow$) By hypothesis there exist nonzero $a, b \in A$ prime to $f$ such that $b/a$ is positive, $a \equiv b \bmod f$ and $aI = bJ$. We have

$$ba^{-1}I^{-1} \cdot W = J^{-1} \cdot W, \quad ba^{-1}x - x \in a^{-1} \cdot W \cap J^{-1}f^{-1} \cdot W \subseteq J^{-1} \cdot W,$$

whence the result.

(I)($\Leftarrow$) By hypothesis there exists some positive $y \in \mathbf{K}$ such that

$$yI^{-1} \cdot W = J^{-1} \cdot W, \quad yx - x \in J^{-1} \cdot W.$$

Necessarily $yJ = I$, and moreover, because $I$ and $J$ are prime to $f$, we can write $y = b/a$ where $0 \neq a, b \in A$ are prime to $f$. We have

$$(b - a)x \in (I^{-1} + J^{-1}) \cdot W \cap f^{-1} \cdot W = W,$$

hence $a \equiv b \bmod f$, and the result follows.

(II) Replacing $x'+W'$ by a homothetic torsion $A$-lattice translate, we may assume that $W' = J^{-1} \cdot W$ for some integral $A$-ideal $J$ prime to $f$; replacing $W$ by $J^{-1} \cdot W$, we may assume that $W = W'$. By hypothesis both $x + W$ and $x' + W$ generate the free rank one $(A/f)$-module $f^{-1} \cdot W/W$, and hence we can find positive $a \in A$ such that $ax \equiv x' \bmod W$. Necessarily $a$ is prime to $f$. Then $x' + W = ax + W \sim x + a^{-1}W = x + (a^{-1}) \cdot W$. ◇

**Proposition 3.4.3** *There exists a unique (left) action of the group $G$ on $\Xi$ such that for all integral $A$-ideals $f$, integral $A$-ideals $I$ prime to $f$, and $A$-lattice translates $x + W$ of order $f$,*

$$\sigma[x + W] = [x + I^{-1} \cdot W]$$

*for all $\sigma \in G$ with image in $G_f$ equal to the $\sim_f$-equivalence class of $I$. (Hereafter $\Xi$ is considered to be equipped with the action of $G$ so defined.)*

**Proof** This follows in a straightforward way from Proposition 3.3.1 and Proposition 3.4.2. ◇

**Proposition 3.4.4** *Let $f$ be an integral $A$-ideal. The set $\Xi(f)$ is $G$-stable. The set $\Xi^\times(f)$ consists of a single $G$-orbit. The isotropy subgroup of each element of $\Xi^\times(f)$ is $\ker(G \to G_f)$. The map $Y_f : \Xi \to \Xi$ is $G$-equivariant.*

**Proof** This follows in a straightforward way from Proposition 3.3.1 and Proposition 3.4.2. ◇

**Remark** In class field theory the group $G$ is identified with the Galois group of a certain infinite abelian extension of $k$. In the archimedean (resp. nonarchimedean) case, this extension can be obtained explicitly by adjoining to $k$ all roots of unity (resp. all torsion points of all sign-normalized rank one elliptic $A$-modules). The set $\Xi$ turns out to be $G$-equivariantly isomorphic in the archimedean (resp. nonarchimedean) case to the set of roots of unity (resp. the disjoint union, extended over the set of sign-normalized rank one elliptic $A$-modules $\rho$, of the set of torsion points of $\rho$). For an overview of the theory of sign-normalized rank one elliptic $A$-modules see Hayes [8]. For an overview of related function field arithmetic see Goss [7].

## 3.5 The sign group $G_\infty$

For each integral $A$-ideal $f$, let $D_f$ be the subgroup of $G_f$ consisting of the $\sim_f$-equivalence classes of ideals of the form $(a)$ for some $a \in A$ such that $a \equiv 1 \bmod f$. We define

$$G_\infty := \varprojlim D_f \subset G.$$

We call $G_\infty$ the *sign group*.

**Remark** Identifying $G$ with the Galois group of an abelian extension of $k$ via class field theory, the subgroup $G_\infty$ may be interpreted as the decomposition group of the valuation of $k$ inherited from $\mathbf{K}$.

**Proposition 3.5.1** *There exists a unique homomorphism*

$$\operatorname{sgn} : G_\infty \to \mathbf{K}^\times$$

*mapping $G_\infty$ isomorphically to the group of roots of unity of $\mathbf{K}$ such that*

$$\gamma[x + W] = [(\operatorname{sgn} \gamma)^{-1}(x + W)] \tag{1}$$

*for all $\gamma \in G_\infty$ and torsion $A$-lattice translates $x + W$.*

**Proof** For each integral $A$-ideal $f$, let $E_f$ be the subgroup of $A^\times$ consisting of $a \in A^\times$ such that $a \equiv 1 \bmod f$. Then for all but finitely many integral $A$-ideals $f$, one has $E_f = \{1\}$. (In fact, if $f$ is not the unit ideal and $E_f \neq \{1\}$, then $(\mathbf{K}, A, \operatorname{sgn}) = (\mathbf{R}, \mathbf{Z}, x \mapsto x/|x|)$ and $f = (2)$.) For all integral $A$-ideals $f$ there exists a unique homomorphism $\phi_f : D_f \to \mathbf{K}^\times/E_f$ mapping $D_f$ isomorphically to the group of roots of unity of $\mathbf{K}^\times$ modulo $E_f$ under which an element of $D_f$ represented by an ideal of the form $(a)$ for some $a \in A$ such that $a \equiv 1 \bmod f$ maps to the $E_f$-coset containing $\operatorname{sgn} a$. The system $\{\phi_f\}$ is compatible and induces a homomorphism $\phi : G_\infty \to \mathbf{K}^\times$ mapping $G_\infty$ isomorphically to the group of roots of unity of $\mathbf{K}$.

We claim that $\phi$ has property (1). Fix an integral $A$-ideal $f$, a torsion $A$-lattice translate $x + W$ of order $f$ and an element $\gamma \in G_\infty$. Choose an integral $A$-ideal $I$ prime to $f$ belonging to the $\sim_f$-equivalence class to which $\gamma$ gives rise in $G_f$. Then $I = (a)$ for some $a \in A$ such that $a \equiv 1 \bmod f$, and

$$\begin{aligned}
\gamma[x + W] &= [x + I^{-1} \cdot W] \\
&= [(\operatorname{sgn} a)^{-1} a (x + I^{-1} \cdot W)] \\
&= [(\operatorname{sgn} a)^{-1}(ax + W)] \\
&= [(\operatorname{sgn} a)^{-1}(x + W)] = [\phi(\gamma)^{-1}(x + W)].
\end{aligned}$$

The claim is proved. Thus we have established the existence of a homomorphism $G_\infty \to \mathbf{K}^\times$ with the desired properties; uniqueness follows by Proposition 3.4.4. ⋄

Fix a generator $\gamma_0 \in G_\infty$ arbitrarily and let $m$ denote the order of $G_\infty$. Given an abelian group $M$ equipped with an action of $G_\infty$ and an integer $i$, we define the $i^{th}$ *sign-homology module* and the $(1-i)^{th}$ *sign-cohomology module* of $M$ to be

$$H_i\left(\mathcal{T}\left(M/\mathbf{Z}[G_\infty], \begin{bmatrix} \sum_{i=0}^{m-1} \gamma_0^i \\ 1 - \gamma_0 \end{bmatrix}\right)\right).$$

As explained in §2, the sign-(co)homology of $M$ can be identified with the Tate (co)homology of $G_\infty$ with coefficients in $M$.

## 3.6 The module $U^{(\nu)}$

Let $R$ be a commutative ring with unit. Let $\mathcal{A}$ be the free $R$-module generated by $\Xi$, and let the action of $G$ on $\Xi$ be extended to $\mathcal{A}$ in $R$-linear fashion. Fix a family

$$\nu = \{\nu_f\}$$

of elements of $R$ indexed by the integral $A$-ideals such that

$$\nu_1 = 1$$

and

$$\nu_{fg} = \nu_f \nu_g$$

for all integral $A$-ideals $f$ and $g$.

The $R$-module $U^{(\nu)}$ is defined to be the quotient of $\mathcal{A}$ by the $R$-submodule generated by the family of elements of the form

$$\nu_p \xi - \sum_{\eta : Y_p \eta = \xi} \eta$$

where $p$ is a maximal $A$-ideal, $\xi \in \Xi$, and the sum is extended over those $\eta \in \Xi$ such that $Y_p \eta = \xi$. By Proposition 3.4.4, $\ker\left(\mathcal{A} \to U^{(\nu)}\right)$ is $G$-stable, and hence the action of $G$ on $\mathcal{A}$ descends to $U^{(\nu)}$. The multiplicativity of the system $\{\nu_f\}$, along with Proposition 3.4.4, implies that $\ker\left(\mathcal{A} \to U^{(\nu)}\right)$ contains every element of the form

$$\nu_f \xi - \sum_{\eta: Y_f \eta = \xi} \eta$$

where $f$ is an integral $A$-ideal, $\xi \in \Xi$, and the sum is extended over those $\eta \in \Xi$ such that $Y_f \eta = \xi$.

Let $f$ be an integral $A$-ideal. We define $\mathcal{A}(f)$ to be the $R$-submodule of $\mathcal{A}$ generated by $\Xi(f)$. The $R$-module $\mathcal{A}(f)$ is finitely generated and free, and clearly

$$\mathcal{A} = \bigcup_f \mathcal{A}(f).$$

By Proposition 3.4.4, the $R$-module $\mathcal{A}(f)$ is a $G$-stable $R$-submodule of $\mathcal{A}$, and moreover the action of $G$ on $\mathcal{A}(f)$ factors through an action of $G_f$. We define $U^{(\nu)}(f)$ to be the quotient of $\mathcal{A}(f)$ by the $R$-submodule generated by all elements of the form

$$\nu_p \xi - \sum_{\eta: Y_p \eta = \xi} \eta$$

where $p$ is an maximal $A$-ideal dividing $f$, $\xi \in \Xi(f/p)$, and the sum is extended over those $\eta \in \Xi$ such that $Y_p \eta = \xi$. By Proposition 3.4.4, $\ker\left(\mathcal{A}(f) \to U^{(\nu)}(f)\right)$ is $G$-stable, and hence the action of $G$ descends to $U^{(\nu)}(f)$. Note that the action of $G$ on $U^{(\nu)}(f)$ factors through an action of $G_f$. Clearly

$$U^{(\nu)} = \varinjlim U^{(\nu)}(f).$$

The multiplicativity of the system $\{\nu_f\}$ implies that $\ker\left(\mathcal{A}(f) \to U^{(\nu)}(f)\right)$ contains all elements of the form

$$\nu_g \xi - \sum_{\eta: Y_g \eta = \xi} \eta$$

where $g$ is any integral $A$-ideal dividing $f$, $\xi \in \Xi(f/g)$, and the sum is extended over $\eta \in \Xi$ such that $Y_g \eta = \xi$.

# 4 The structure of $U^{(\nu)}$ and its sign-homology

## 4.1 A partition of $\Xi$

**Lemma 4.1.1** *Let $f$ be an integral $A$-ideal. Let $p$ be a maximal $A$-ideal dividing $f$. Write $f = cp^n$ where $n$ is a positive integer and $c$ is an integral $A$-ideal prime to $f$. Let $\phi \in G$ be an element projecting to the $\sim_c$-equivalence class of $p$ in $G_c$. Let $S$ be a set of elements of $G$ mapping bijectively to $\ker\left(G_f \to G_{f/p}\right)$ under projection to $G_f$. Then for each $\xi \in \Xi^\times(f)$ one has*

$$\sum_{\eta: Y_p\eta = Y_p\xi} \eta = \left(\sum_{\sigma \in S} \sigma\xi\right) + \begin{cases} \phi^{-1} Y_p\xi & \text{if } n = 1 \\ 0 & \text{if } n > 1 \end{cases} \tag{2}$$

*where the sum on the left is extended over $\eta \in \Xi$ such that $Y_p\eta = Y_p\xi$.*

**Proof** Let $x + W$ be an $A$-lattice translate of order $f$ such that $\xi = [x + W]$. Let $T$ be a set of positive elements of $A$ prime to $f$ mapping bijectively to $\ker\left((A/f)^\times \to (A/(f/p))^\times\right)$ under reduction modulo $f$. Let $b$ be a positive element of $A$ such that $b \equiv 0 \bmod p$ and $b \equiv 1 \bmod c$. Then the sum

$$\left(\sum_{a \in T}[ax + W]\right) + \begin{cases} [bx + W] & \text{if } n = 1 \\ 0 & \text{if } n > 1 \end{cases} \tag{3}$$

equals the left side of (2). Put $J := p^{-1}(b)$. Then $J$ is an integral $A$-ideal prime to $c$ such that $Jp$ is $\sim_c$-equivalent to the unit ideal, $\{(a) \mid a \in T\}$ is a set of integral $A$-ideals prime to $f$ mapping bijectively to $\ker\left(G_f \to G_{f/p}\right)$, and

$$\left(\sum_{a \in T}[x + (a)^{-1} \cdot W]\right) + \begin{cases} [x + J^{-1}p^{-1} \cdot W] & \text{if } n = 1 \\ 0 & \text{if } n > 1 \end{cases} \tag{4}$$

equals the right side of (2). But for each $a \in T$ one has

$$ax + W \sim x + a^{-1}W = x + (a)^{-1} \cdot W.$$

Further, in the case $n = 1$ one has

$$bx + w \sim x + b^{-1}W = x + J^{-1}p^{-1} \cdot W.$$

Thus the sums (3) and (4) are equal term by term. ◇

For each integral $A$-ideal $f$, put $\Xi(f^\infty) := \bigcup_{N=1}^\infty \Xi(f^N)$.

**Lemma 4.1.2** *There exists a partition*

$$\Xi = \coprod_{k=0}^\infty \Xi_k$$

*with the following properties:*

(I) For all positive integers $k$, integral $A$-ideals $f$, and $\xi \in \Xi_k \cap \Xi^\times(f)$, there exists a maximal $A$-ideal $p$ dividing $f$ such that for all $\eta \in \Xi$, if $Y_p \eta = Y_p \xi$ and $\xi \neq \eta$, then either $\eta \in \Xi_{k-1}$ or $\eta \in \Xi(f/p)$.

(II) For all integral $A$-ideals $f$, with the exception in the archimedean case of those $f$ exactly divisible by 2, the sign group $G_\infty$ stabilizes and acts freely on the set $\Xi_0 \cap (\Xi(f^\infty) \setminus \Xi(f))$.

(III) For each integral $A$-ideal $f$ one has $|\Xi_0 \cap \Xi(f)| = |G_f|$.

(Hereafter we will assume such a partition of $\Xi$ to be fixed.)

**Proof** For each integral $A$-ideal $f$, we select a subset $S(f) \subseteq G$ with the following properties:

- $1 \in S(f)$.

- The natural map $S_f \to G_f$ is bijective.

- If the natural map $G_\infty \to G_f$ is injective, then the set $S(f)$ is stable under the action of $G_\infty$.

By Proposition 3.4.1 and Proposition 3.4.4 it follows that for each integral $A$-ideal $f$ and $\xi \in \Xi^\times(f)$, there exist unique

$$\sigma \in S(f/\sqrt{f})$$

and, for each maximal $A$-ideal $p$ dividing $f$, unique

$$\tau_p \in S(f) \cap \ker\left(G \to G_{f/p}\right)$$

such that

$$\xi = \sigma \left(\prod_p \tau_p\right)[1 + f];$$

in this situation we declare $\xi \in \Xi_k$, where $k$ is the number of maximal $A$-ideals $p$ dividing $f$ such that $\tau_p = 1$. Property (I) of the partition so defined is easily verified with the help of Lemma 4.1.1. We claim that for all integral $A$-ideals $g$ such that $\sqrt{g}$ divides $\sqrt{f}$ but $g$ does not divide $f$, the natural map $G_\infty \to G_{g/\sqrt{g}}$ is injective. The proof of the claim is straightforward (if somewhat tedious) and so we omit it. The claim proves (II). For each integral $A$-ideal $f$, one has

$$|\Xi_0 \cap \Xi^\times(f)| = \left|G_{f/\sqrt{f}}\right| \cdot \prod_{p|f} \left(\left|\ker\left(G_f \to G_{f/p}\right)\right| - 1\right)$$

where the product is extended over the maximal $A$-ideals $p$ dividing $f$, whence follows property (III) of the partition by Möbius inversion. ⋄

## 4.2  A $\Lambda$-basis for $\mathcal{A}$

Let $\Lambda$ be the polynomial ring over $R$ generated by a collection $\{X_p\}$ of independent variables indexed by the maximal $A$-ideals. For each integral $A$-ideal $f$ put

$$X_f := \prod_i X_{p_i}^{n_i}$$

where

$$f = \prod_i p_i^{n_i}$$

is the prime factorization of $f$. Then the collection of monomials $\{X_f\}$ indexed by the integral $A$-ideals is an $R$-basis of $\Lambda$. We equip $\mathcal{A}$ with the unique structure of $\Lambda$-module extending the $R$-module structure in such a way that

$$X_f \xi = \sum_{\eta : Y_f \eta = \xi} \eta$$

for all $\xi \in \Xi$ and integral $A$-ideals $f$, where the sum is extended over $\eta \in \Xi$ such that $Y_f \eta = \xi$. By Proposition 3.4.4, the action of $G$ on $\mathcal{A}$ is $\Lambda$-linear.

For each integral $A$-ideal $f$, let $\mathcal{A}(f^\infty)$ be the $R$-span of $\Xi(f^\infty)$, and let $R\left[\{X_p\}_{p|f}\right]$ be the $R$-subalgebra of $\Lambda$ generated by the variables $X_p$ where $p$ is a maximal $A$-ideal dividing $f$. Note that $\mathcal{A}(f^\infty)$ is a $G$-stable $R\left[\{X_p\}_{p|f}\right]$-submodule of $\mathcal{A}$.

**Theorem 4.2.1** *Let $f$ be an integral $A$-ideal.*

*(I) The elements of $\mathcal{A}(f)$ of the form $X_g \xi$ with $g$ an integral $A$-ideal dividing $f$ and $\xi \in \Xi_0 \cap \Xi(f/g)$ constitute an $R$-basis.*

(II) *The elements of $\mathcal{A}(f^\infty)$ of the form $X_g \xi$ with $g$ a integral $A$-ideal such that $\sqrt{g}$ divides $\sqrt{f}$ and $\xi \in \Xi_0 \cap \Xi(f^\infty)$ constitute an $R$-basis.*

(III) *As an $R\left[\{X_p\}_{p|f}\right]$-module $\mathcal{A}(f^\infty)$ is free and the set $\Xi_0 \cap \Xi(f^\infty)$ is a basis.*

(IV) *As a $\Lambda$-module $\mathcal{A}$ is free and the set $\Xi_0$ is a basis.*

**Proof** Clearly (I)⇒(II)⇒(III)⇒(IV). It will be enough to prove (I). In turn, it will be enough to show that the family of elements of $\mathcal{A}(f)$ in question spans $\mathcal{A}(f)$ over $R$, because that family has cardinality

$$\sum_{g|f} |G_g| = \sum_{g|f} |\Xi^\times(g)| = |\Xi(f)|$$

by property (I) of the partition $\Xi = \coprod_{k=0}^\infty \Xi_k$. In turn, it will be enough to prove that

$$\mathcal{A}_k \cap \mathcal{A}^\times(f) \subseteq \left(\mathcal{A}_{k-1} \cap \mathcal{A}^\times(f)\right) + \sum_{p|f} \mathcal{A}(f/p) + \sum_{p|f} X_p \mathcal{A}(f/p) \tag{5}$$

where $k$ is a positive integer, $\mathcal{A}_k$ is the $R$-span of $\Xi_k$, $\mathcal{A}^\times(f)$ is the $R$-span of $\Xi^\times(f)$, and the sums are extended over the maximal $A$-ideals $p$ dividing $f$. But for each $\xi \in \Xi^\times(f) \cap \Xi_k$ there exists by property (I) of the partition $\Xi = \coprod_{k=0}^\infty \Xi_k$ a maximal $A$-ideal $p$ such that

$$\xi - X_p Y_p \xi \in \mathcal{A}^\times \cap \mathcal{A}_{k-1}(f) + \mathcal{A}(f/p),$$

and hence (5) does indeed hold. ◇

## 4.3 An $R$-basis for $U^{(\nu)}$

For each integral $A$-ideal $f$, put

$$U^{(\nu)}(f^\infty) := \lim_{N \to \infty} U^{(\nu)}(f^N) = \frac{\mathcal{A}(f^\infty)}{\sum_{p|f}(X_p - \nu_p)\mathcal{A}(f^\infty)}$$

where the sum is extended over maximal $A$-ideals $p$ dividing $f$.

**Theorem 4.3.1** *Let $f$ be an integral $A$-ideal.*

(I) *The $R$-module $U^{(\nu)}(f)$ is free and the set $\Xi_0 \cap \Xi(f)$ gives rise to an $R$-basis.*

(II) The $R$-module $U^{(\nu)}(f^\infty)$ is free and the set $\Xi_0 \cap \Xi(f^\infty)$ gives rise to an $R$-basis.

(III) With the exception in the archimedean case of $f$ exactly divisible by 2, the natural map $U^{(\nu)}(f) \to U^{(\nu)}(f^\infty)$ induces an isomorphism in sign-homology.

(IV) The $R$-module $U^{(\nu)}$ is free and the set $\Xi_0$ gives rise to an $R$-basis.

**Proof** (IV) Clearly

$$R = \frac{\Lambda}{\sum_p (X_p - \nu_p)\Lambda}, \quad U^{(\nu)} = \frac{\mathcal{A}}{\sum_p (X_p - \nu_p)\mathcal{A}}.$$

By Theorem 4.2.1, the set $\Xi_0$ is a $\Lambda$-basis of $\mathcal{A}$, and therefore gives rise to an $R$-basis of $U^{(\nu)}$.

(I) Let $\mathcal{A}_0$ be the $R$-span of $\Xi_0 \cap \Xi$. It is enough to show that the natural map $\mathcal{A}_0 \cap \mathcal{A}(f) \to U^{(\nu)}(f)$ is bijective, and injectivity is clear by what we have proved so far. By Theorem 4.2.1 we have

$$\mathcal{A}(f) = (\mathcal{A}_0 \cap \mathcal{A}(f)) \bigoplus \left( \sum_{1 \neq g | f} X_g (\mathcal{A}_0 \cap \mathcal{A}(f/g)) \right)$$

and hence

$$\mathcal{A}(f) = (\mathcal{A}_0 \cap \mathcal{A}(f)) \bigoplus \left( \sum_{1 \neq g | f} (X_g - \nu_g)(\mathcal{A}_0 \cap \mathcal{A}(f/g)) \right),$$

where the sums are extended over nonunit integral $A$-ideals $g$ dividing $f$, whence follows the surjectivity of map in question.

(II) This is a trivial consequence of (I).

(III) By Lemma 4.1.2 and what we have already proved, the set $\Xi_0 \cap (\Xi(f^\infty) \setminus \Xi(f))$ gives rise to an $R$-basis for the quotient $U^{(\nu)}(f^\infty)/U^{(\nu)}(f)$ that is stabilized by $G_\infty$ and on which $G_\infty$ acts freely. Consequently the sign-homology of the quotient $U^{(\nu)}(f^\infty)/U^{(\nu)}(f)$ vanishes. ⋄

## 4.4 The sign-homology of $U^{(\nu)}$

Let $\mathcal{A}'$ be the $R$-submodule of $\mathcal{A}$ generated by $\Xi \setminus \Xi(1)$. Let $\Lambda[G]$ be the group ring of $G$ with coefficients in $\Lambda$. Note that $\mathcal{A}$ is a $\Lambda[G]$-module and that $\mathcal{A}'$ is a $\Lambda[G]$-submodule. Let $R[G] \subset \Lambda[G]$ be the $R$-subalgebra generated by $G$ and let $R[G_\infty] \subseteq \Lambda[G]$ be the $R$-subalgebra generated by $G_\infty$. Let $\Xi^\dagger \subset \Xi$ be the union of all $G_\infty$-orbits of cardinality $|G_\infty|$. Let $\mathcal{A}^\dagger$ be the $R$-span of $\Xi^\dagger$. Then $\mathcal{A}^\dagger$ is a free $R[G_\infty]$-module. Note that $\mathcal{A}^\dagger$ is a $\Lambda[G]$-submodule of $\mathcal{A}$. Note that $\mathcal{A}'$ is a $\Lambda[G]$-submodule of $\mathcal{A}$ containing $\mathcal{A}^\dagger$. Of course the only case in which $\mathcal{A}' \neq \mathcal{A}^\dagger$ is the archimedean case.

**Proposition 4.4.1** *In the archimedean case, $X_p$ annihilates $\mathcal{A}'/\mathcal{A}^\dagger$ for all primes $p$.*

**Proof** One has $X_2[1/2 + \mathbf{Z}] = [1/4 + \mathbf{Z}] + [3/4 + \mathbf{Z}]$. The remaining cases are similarly trivial. $\diamond$

Fix a linear ordering of the set of maximal $A$-ideals arbitrarily.

**Theorem 4.4.2** *Assume either that we are in the nonarchimedean case or that $\nu_f \in R^\times$ for all integral $A$-ideals $f$. Let $m$ be the order of $G_\infty$ and let $\gamma_0 \in G_\infty$ be a generator. Then the directed family of graded $R[G]$-modules underlying the directed family*

$$\left\{ H_* \left( \mathcal{KT}^{tot} \left( (\mathcal{A}/\mathcal{A}')/\Lambda[G], \{\nu_p - X_p\}_{p|f}, \begin{bmatrix} \sum_{i=0}^{m-1} \gamma_0^i \\ 1 - \gamma_0 \end{bmatrix} \right) \right) \right\} \qquad (6)$$

*of graded $\Lambda[G]$-modules indexed by squarefree integral $A$-ideals $f$ is isomorphic to the directed family*

$$\left\{ H_* \left( \mathcal{T} \left( U^{(\nu)}(f^\infty)/R[G], \begin{bmatrix} \sum_{i=0}^{m-1} \gamma_0^i \\ 1 - \gamma_0 \end{bmatrix} \right) \right) \right\} \qquad (7)$$

*of graded $R[G]$-modules indexed by squarefree integral $A$-ideals $f$. (An explicit isomorphism is given in the proof.)*

**Proof** Let $f$ be a squarefree integral $A$-ideal. Let $\Lambda(f)[G]$ be the group ring of $G$ with coefficients in the $R$-subalgebra $\Lambda(f) \subseteq \Lambda$ generated by the variables $X_p$ for $p$ ranging over maximal $A$-ideals dividing $f$. Consider the following chain complexes of $\Lambda(f)[G]$-modules.

$$\mathcal{KT}^{tot} \left( (\mathcal{A}/\mathcal{A}')/\Lambda(f)[G], \{\nu_p - X_p\}_{p|f}, \begin{bmatrix} \sum_{i=0}^{m-1} \gamma_0^i \\ 1 - \gamma_0 \end{bmatrix} \right) \qquad (8)$$

$$\mathcal{KT}^{tot}\left(\mathcal{A}(f^\infty)/\Lambda(f)[G], \{\nu_p - X_p\}_{p|f}, \begin{bmatrix} \sum_{i=0}^{m-1} \gamma_0^i \\ 1 - \gamma_0 \end{bmatrix}\right) \tag{9}$$

$$\overline{\mathcal{T}}\left(\mathcal{A}(f^\infty)/\Lambda(f)[G], \{\nu_p - X_p\}_{p|f}, \begin{bmatrix} \sum_{i=0}^{m-1} \gamma_0^i \\ 1 - \gamma_0 \end{bmatrix}\right) \tag{10}$$

The chain complex (8) is naturally a quotient of (9) because

$$\mathcal{A}/\mathcal{A}' = \mathcal{A}(f^\infty)/(\mathcal{A}(f^\infty) \cap \mathcal{A}').$$

The chain complex (10) is naturally a quotient of (9) because they are companions of the double complex

$$\mathcal{KT}\left(\mathcal{A}(f^\infty)/\Lambda(f)[G], \{\nu_p - X_p\}_{p|f}, \begin{bmatrix} \sum_{i=0}^{m-1} \gamma_0^i \\ 1 - \gamma_0 \end{bmatrix}\right). \tag{11}$$

We claim that both quotient maps induce homology isomorphisms; the claim granted, the isomorphism from the homology of (8) to the homology of (10) provided by the claim induces (for variable $f$) the desired isomorphism from the directed family of graded $R[G]$-modules underlying (6) to the directed family (7).

We turn to the proof of the claim. Clearly the chain complex

$$\mathcal{T}\left(\mathcal{A}(f^\infty) \cap \mathcal{A}^\dagger/\Lambda(f)[G], \begin{bmatrix} \sum_{i=0}^{m-1} \gamma_0^i \\ 1 - \gamma_0 \end{bmatrix}\right)$$

is acyclic, and hence by Proposition 2.4.1 the chain complex

$$\mathcal{KT}^{tot}\left(\mathcal{A}(f^\infty) \cap \mathcal{A}^\dagger/\Lambda(f)[G], \{\nu_p - X_p\}_{p|f}, \begin{bmatrix} \sum_{i=0}^{m-1} \gamma_0^i \\ 1 - \gamma_0 \end{bmatrix}\right) \tag{12}$$

is acyclic. The chain complex

$$\mathcal{K}\left(\left(\frac{\mathcal{A}(f^\infty) \cap \mathcal{A}'}{\mathcal{A}(f^\infty) \cap \mathcal{A}^\dagger}\right)/\Lambda(f)[G], \{\nu_p - X_p\}_{p|f}\right)$$

vanishes in the nonarchimedean case and is acyclic in the archimedean case by Proposition 2.4.1 and Proposition 4.4.1. Therefore the chain complex

$$\mathcal{KT}^{tot}\left(\left(\frac{\mathcal{A}(f^\infty) \cap \mathcal{A}'}{\mathcal{A}(f^\infty) \cap \mathcal{A}^\dagger}\right)/\Lambda(f)[G], \{\nu_p - X_p\}_{p|f}, \begin{bmatrix} \sum_{i=0}^{m-1} \gamma_0^i \\ 1 - \gamma_0 \end{bmatrix}\right) \tag{13}$$

is acyclic by Proposition 2.4.1. The acyclicity of the chain complexes (12) and (13) implies that the quotient map from (9) to (8) is a homology isomorphism. Now $\mathcal{A}(f^\infty)$ is by Theorem 4.2.1 a free $\Lambda(f)[G]$-module. It follows by Proposition 2.4.1 that the complex

$$\mathcal{K}\left(\mathcal{A}(f^\infty)/\Lambda(f)[G], \{\nu_p - X_p\}_{p|f}\right)$$

is acyclic in positive degree, and hence by Proposition 2.4.1 the quotient map (9) to (10) is indeed a homology isomorphism. The claim is proved. ⋄

Let $\tilde{\Lambda}$ be the ring obtained from $\Lambda$ by inverting the variables $X_p$ for $p$ ranging over maximal $A$-ideals, and let $\Gamma \subset \tilde{\Lambda}^\times$ be the subgroup generated by those variables. Then $\tilde{\Lambda}$ is the group ring of $\Gamma$ over $R$, and $\Gamma$ is a free abelian group for which the family of elements of the form $X_p$ for $p$ a maximal $A$-ideal constitute a basis. Let $\tilde{\Lambda}[G]$ be the group ring of $G$ with coefficients in $\tilde{\Lambda}$, and let $\Gamma G \subset \tilde{\Lambda}[G]^\times$ be the subgroup generated by $\Gamma$ and $G$. Then the natural map $\Gamma \times G \to \Gamma G$ is an isomorphism and $\tilde{\Lambda}[G]$ may be viewed as the group ring of $\Gamma G$ with coefficients in $R$. Let $\Pi \subseteq \Gamma G$ be the subgroup generated by monomials of the form $X_p\gamma$ with $p$ a maximal $A$-ideal and $\gamma$ an element of $G$ such that $p$ and $\gamma$ have a common image in the narrow ideal class group $G_1$.

**Proposition 4.4.3** *For every maximal $A$-ideal $p$, the action of $X_p$ on $\mathcal{A}/\mathcal{A}'$ is invertible, and hence $\mathcal{A}/\mathcal{A}'$ can be viewed as a $\Gamma G$-module; as such, $\mathcal{A}/\mathcal{A}'$ is isomorphic to $\mathrm{Ind}_\Pi^{\Gamma G} \mathrm{Inf}_{\{1\}}^\Pi R$.*

**Proof** For any $A$-lattice $W$, maximal $A$-ideals $p$ and $q$, and $\gamma \in G$ such that $\gamma$ and $q$ have a common image in $G_1$, one has

$$\gamma X_p[W] \;=\; [q^{-1}p \cdot W] + \sum_{0 \neq y \in q^{-1}p \cdot W / p \cdot W} [y + p \cdot W]$$

$$=\; [q^{-1}p \cdot W] \bmod \mathcal{A}'.$$

The result follows straightforwardly from this identity. ⋄

# 5 The universal ordinary distribution

For the remainder of the paper we specialize the preceding theory as follows. We take the coefficient ring $R$ to be $\mathbf{Z}$ and we take $\nu_f = 1$ for all integral $A$-ideals $f$. Then $\mathcal{A}$ becomes the free abelian group generated by the set $\Xi$ and $U$ becomes the quotient of $\mathcal{A}$ by the subgroup generated by all elements of the form $\xi - \sum_{\eta:Y_p\eta=\xi} \eta$ with $\xi \in \Xi$ and $p$ a maximal $A$-ideal. We now write simply $U$ instead of $U^{(\nu)}$ and $U(f)$ instead of $U^{(\nu)}(f)$. We call $U$ the *universal ordinary distribution*.

## 5.1 Comparison of $U(f)$ and $U'(f)$

For each locally constant homomorphism $\chi : G \to \mathbf{C}^\times$, there exists an integral $A$-ideal $c$ such that for all integral $A$-ideals $f$, the homomorphism $\chi$ factors through $G_f$ if and only if $c$ divides $f$; the integral $A$-ideal $c$ is called the *conductor* of $\chi$. Given a locally constant homomorphism $\chi : G \to \mathbf{C}^\times$ of conductor $c$ and a maximal $A$-ideal $p$, if $p$ does not divide $c$, let $\chi(p)$ denote the value of $\chi$ at any $\sigma \in G$ projecting to the $\sim_c$-equivalence class of $p$ in $G_c$, and otherwise, if $p$ does divide $c$, put $\chi(p) := 0$.

**Lemma 5.1.1** *There exists a unique homomorphism $u : \mathcal{A} \to \mathbf{Q}$ such that for all integral $A$-ideals $f$ and locally constant homomorphisms $\chi : G \to \mathbf{C}^\times$ of conductor dividing $f$ one has*

$$\int_G u(\gamma[1+f])\chi(\gamma)d\mu(\gamma) := \frac{1}{|G_f|} \prod_{p|f} (1 - \chi(p))$$

*where $\mu$ is Haar probability measure on $G$ and the product is extended over the maximal $A$-ideals $p$ dividing $f$. Necessarily $u$ factors through the universal ordinary distribution $U$.*

**Proof** Existence and uniqueness of $u$ are clear in view of Proposition 3.4.4. Fix a maximal $A$-ideal $p$ and consider the unique homomorphism $v : \mathcal{A} \to \mathbf{Q}$ such that $v(\xi) = u(X_p\xi)$ for all $\xi \in \Xi$. It will be enough to prove that $u = v$. Let $f$ be any integral $A$-ideal divisible by $p$. Write $f = cp^n$ with $n \geq 1$ and $c$ prime to $p$. Select a subset $S \subset G$ mapping bijectively to $\ker\left(G_f \to G_{f/p}\right)$, and select $\phi \in G$ such that $\phi$ and $p$ have a common image in the generalized ideal class group $G_c$. For all locally constant homomorphisms $\chi : G \to \mathbf{C}^\times$ of conductor dividing $f/p$, by Lemma 4.1.1 and the definition of $u$, one has

$$\int_G v(\gamma[1+f/p])\chi(\gamma)d\mu(\gamma)$$

$$= \int_G \left(\sum_{\sigma \in S} u(\sigma\gamma[1+f]) + \begin{cases} u(\phi^{-1}\gamma[1+f/p]) & \text{if } n=1 \\ 0 & \text{if } n>1 \end{cases}\right)\chi(\gamma)d\mu(\gamma)$$

$$= \frac{1}{|G_{f/p}|}\prod_{q|(f/p)}(1-\chi(q))$$

where the product is extended over the maximal $A$-ideals $q$ dividing $f$. Therefore $u = v$, and we are done. ⋄

For each integral $A$-ideal $f$, let $U'(f)$ be the $\mathbf{Z}[G_f]$-submodule of $\mathbf{Q}[G_f]$ generated by elements of the form

$$\sum_{\gamma \in G_f} u(\tilde{\gamma}[1+g])\gamma^{-1} \in \mathbf{Q}[G_f]$$

where $g$ is any integral $A$-ideal dividing $f$ and $\tilde{\gamma} \in G$ is any lifting of $\gamma \in G_f$. It is easy to check that $U'(f)$ is a free abelian group of rank $|G_f|$. In the archimedean case $U'(f)$ coincides with the module denoted by $U$ in Sinnott's paper [15], and in the nonarchimedean case with the module denoted by $U$ in Yin's paper [16]. In all cases, because the homomorphism $u$ factors through the universal ordinary distribution $U$, it follows formally that $U'(f)$ is $G$-equivariantly a quotient of $U(f)$ and hence $G$-equivariantly isomorphic to $U(f)$ because the underlying abelian groups are free of the same rank, namely $|G_f|$.

## 5.2 Proof of Yin's conjecture

Assume now that we are in the nonarchimedean case. The following was conjectured L. S. Yin [16].

**Theorem 5.2.1** *Let $f$ be a nonunit integral $A$-ideal. Let $r$ be the number of distinct maximal $A$-ideals dividing $f$. Let $H_f \subseteq G_1$ be the subgroup generated by the narrow ideal classes of the maximal $A$-ideals dividing $f$. Let $w$ be the number of roots of unity of $k$. Then the sign-homology of $U(f)$ in each degree a free $(\mathbf{Z}/w\mathbf{Z})[G_1/H_f]$-module of rank $2^{r-1}$.*

**Proof** By Theorem 4.3.1 we may assume that $f$ is squarefree, and it will be enough to show that the sign-homology of $U(f^\infty)$ is in each degree a free $(\mathbf{Z}/w\mathbf{Z})[G_1/H_f]$-module of rank $2^{r-1}$. According to Theorem 4.4.2 we can identify the sign-homology of $U(f^\infty)$ as a graded $G$-module with the graded $G$-module underlying the homology of the chain complex

$$\mathcal{KT}^{tot}\left((\mathcal{A}/\mathcal{A}')/\Lambda[G], \{1-X_p\}_{p|f}, \begin{bmatrix}\sum_{i=0}^{m-1}\gamma_0 \\ 1-\gamma_0\end{bmatrix}\right),$$

where $\gamma_0 \in G_\infty$ is a generator, and in turn, by Proposition 4.4.3, we can identify the homology of the latter complex with the graded $\Lambda[G]$-module underlying the homology of the chain complex

$$\mathcal{KT}^{tot}\left(\mathrm{Ind}_\Pi^{\Gamma G}\mathrm{Inf}_{\{1\}}^\Pi \mathbf{Z}/\mathbf{Z}[\Gamma G], \{1-X_p\}_{p|f}, \begin{bmatrix}\sum_{i=0}^{m-1}\gamma_0 \\ 1-\gamma_0\end{bmatrix}\right).$$

By Proposition 2.5.4 the homology of the complex above is isomorphic as a graded $\Gamma G$-module to

$$\mathrm{Ind}_{\Pi\Gamma(f)G_\infty}^{\Gamma G}\mathrm{Inf}_{\{1\}}^{\Pi\Gamma(f)G_\infty}\hat{H}_*((\Gamma(f)G_\infty)\cap\Pi),$$

where $\Gamma(f)$ is the subgroup of $\Gamma$ generated by $X_p$ for maximal $A$-ideals $p$ dividing $f$. Now the group $(\Gamma(f)G_\infty)\cap\Pi$ is almost free of rank $r$, with torsion subgroup of order $w$, and hence by Proposition 2.5.3,

$$\hat{H}_*((\Gamma(f)G_\infty)\cap\Pi)$$

is in each degree a free $(\mathbf{Z}/w\mathbf{Z})$-module of rank $2^{r-1}$. Finally, one has an isomorphism $G_1/H_f\tilde{\to}\Gamma G/\Pi\Gamma(f)G_\infty$ induced by the inclusion $G\to\Gamma G$. $\diamond$

## 5.3 The archimedean case: the double complex $\mathcal{SK}$

We narrow the focus to the archimedean case. We speak now of positive integers and prime numbers rather than integral and maximal $A$-ideals. We identify $\mathcal{A}$ with the free abelian group on symbols of the form $[a]$ with $a \in \mathbf{Q}/\mathbf{Z}$, and thus identify $U$ with the universal ordinary distribution as defined by Kubert; similarly, for each positive integer $f$, we identify $U(f)$ with Kubert's universal level $f$ ordinary distribution. We also write $U(f^\infty) := \lim_{n\to\infty}U(f^n)$, and we put $\frac{1}{f^\infty}\mathbf{Z}/\mathbf{Z} := \bigcup_{n=1}^\infty \frac{1}{f^n}\mathbf{Z}/\mathbf{Z}$.

By Theorem 4.3.1, we have at our disposal a subset $X_0 \subset \mathbf{Q}/\mathbf{Z}$ such that for all positive integers $f$, the set $X_0 \cap \frac{1}{f}\mathbf{Z}/\mathbf{Z}$ is of cardinality $|G_f|$ and gives rise to a basis of $U(f)$. Thus we recover Kubert's result to the effect that for all positive integers $f$, the natural map $U(f) \to U$ is a split monomorphism with source a free abelian group of rank $|G_f|$. Moreover, by Theorem 4.3.1, for all positive integers $f$ such that $f \not\equiv 2 \bmod 4$, the natural map $U(f) \to U(f^\infty)$ is an isomorphism in sign-homology. In order to recover Kubert's results on the sign-cohomology of $U(f)$ it will be enough to show for each squarefree positive integer $f$ how to use our method to compute the sign-homology of $U(f^\infty)$ and to prove that the natural map $U(f^\infty) \to U$ induces a monomorphism in sign-homology. We will actually do more than this: we will extract a user-friendly description of the sign-homology of $U(f^\infty)$ from the proof of Theorem 4.4.2.

Consider the free abelian group $\mathcal{SK}$ generated by symbols of the form $[a, g, n]$ where $a \in \mathbf{Q}/\mathbf{Z}$, $g$ is a squarefree positive integer and $n$ is an integer. (The notation is meant to call Sinnott and Kubert to mind.) Let $G$ operate on $\mathcal{SK}$ by the rule $\sigma[a, g, n] := [ta, g, n]$, where $t$ is any integer such that for any root of unity $\zeta$ of order equal to the order of $a$, one has $\sigma\zeta = \zeta^t$. Equip $\mathcal{SK}$ with a $G$-stable bigrading $\mathcal{SK} = \bigoplus_m \bigoplus_n \mathcal{SK}_{mn}$ by declaring the symbol $[a, g, n]$ to be of bidegree $(m, n)$, where $m$ is the number of prime factors of $g$. Equip $\mathcal{SK}$ with a $G$-equivariant differential of bidegree $(0, -1)$ by the rule

$$\delta[a, g, n] := (-1)^m([a, g, n-1] + (-1)^n[-a, g, n-1])$$

and a $G$-equivariant differential $\partial$ of bidegree $(-1, 0)$ by the rule

$$\partial[a, g, n] := \sum_{i=1}^{m}(-1)^{i-1}\left([a, g/p_i, n] - \sum_{p_i b = a}[b, g/p_i, n]\right)$$

where $p_1 < \ldots < p_m$ are the primes dividing $g$. Let $\mathcal{SK}' \subset \mathcal{SK}$ be the subgroup generated by symbols $[a, g, n]$ with $a \neq 0$. For each squarefree positive integer $f$, let $\mathcal{SK}(f) \subset \mathcal{SK}$ be the subgroup generated by symbols $[a, g, n]$ where $a \in \frac{1}{f^\infty}\mathbf{Z}/\mathbf{Z}$. Finally, let $\mathcal{N}$ be the subgroup generated by all symbols of the form $[a, g, n]$ with $g \neq 1$, and by all elements of the form $\partial[a, g, n]$ where $g$ is prime. Note that $\mathcal{SK}'$, $\mathcal{SK}(f)$, and $\mathcal{N}$ are bigraded, $G$-, $\partial$- and $\delta$-stable subgroups of $\mathcal{SK}$.

Let $f > 1$ be a squarefree positive integer. On the one hand, the total complex associated to the double complex

$$\frac{\mathcal{SK}(f)}{\mathcal{SK}(f) \cap \mathcal{N}}$$

in an obvious way computes the sign-homology of $U(f^\infty)$. But the double complex somewhat dauntingly denoted

$$\mathcal{KT}\left(\mathcal{A}(f^\infty)/\Lambda(f)[G], \{\nu_p - X_p\}_{p|f}, \begin{bmatrix} \sum_{i=0}^{m-1} \gamma_0^i \\ 1 - \gamma_0 \end{bmatrix}\right).$$

in Theorem 4.4.2 can (as a double complex of $G$-modules) be identified with $\mathcal{SK}(f)$, and what the proof of the theorem says in the present context is that the quotient maps

$$\begin{array}{ccc} \mathcal{SK}(f) & \to & \dfrac{\mathcal{SK}(f)}{\mathcal{SK}(f) \cap \mathcal{SK}'} \\ \downarrow & & \\ \dfrac{\mathcal{SK}(f)}{\mathcal{SK}(f) \cap \mathcal{N}} & & \end{array}$$

induce isomorphisms in homology of associated total complexes. In particular, the total complex associated to the double complex

$$\frac{\mathcal{SK}(f)}{\mathcal{SK}(f) \cap \mathcal{SK}'}$$

also computes the sign-homology of $U(f^\infty)$. But the latter double complex has an extremely simple structure: it is a copy of the double complex

$$\mathcal{KT}\left(\mathbf{Z}/\mathbf{Z}, \{0\}_{p|f}, \begin{bmatrix} 2 \\ 0 \end{bmatrix}\right)$$

which, if employed as in the proof of Proposition 2.5.3, computes the Farrell-Tate homology of the subgroup of $\mathbf{Q}^\times$ generated by $-1$ and the primes $p$ dividing $f$. In particular, as expected, we find that the sign-homology of $U(f^\infty)$ is in each degree a (free) $(\mathbf{Z}/2\mathbf{Z})$-module of rank $2^{r-1}$ where $r$ is the number of prime divisors of $f$.

Passing to the limit over $f$, we can identify the sign-homology of $U$ with the homology of the total complex associated to the double complex $\mathcal{SK}/\mathcal{SK}'$. It is easy to see that the natural map

$$\frac{\mathcal{SK}(f)}{\mathcal{SK}(f) \cap \mathcal{SK}'} \to \frac{\mathcal{SK}}{\mathcal{SK}'}$$

of double complexes is isomorphic to the natural map

$$\mathcal{KT}\left(\mathbf{Z}/\mathbf{Z}, \{0\}_{p|f}, \begin{bmatrix} 2 \\ 0 \end{bmatrix}\right) \to \mathcal{KT}\left(\mathbf{Z}/\mathbf{Z}, \{0\}_{p:\text{any prime}}, \begin{bmatrix} 2 \\ 0 \end{bmatrix}\right)$$

of double complexes. In particular, as expected, it follows that the natural map $U(f^\infty) \to U$ induces an injective map in sign-homology. Finally, it is clear that $G$ operates trivially on the double complex $\mathcal{SK}/\mathcal{SK}'$ and hence $G$ operates trivially, as expected, on the sign-homology of $U$.

**Acknowledgements** I thank W. Sinnott for correspondence about topics related to his work. I thank M. Feshbach for providing me with a reference (namely [1]). I thank G. Mislin for correspondence about his work. I thank the organizers of the conference on Recent Progress in Algebra at KAIST for inviting me and for producing a well run and enjoyable conference. I am especially grateful to my students S. Sinha and P. Das for teaching me a lot about the topics discussed in this paper.

# References

[1] Brown, K. S.: *Cohomology of Groups,* Graduate Texts in Mathematics **87**, Springer-Verlag, New York 1982

[2] Das, P.: *Double coverings of cyclotomic fields arising from algebraic $\Gamma$-monomials.* Ph.D. thesis, University of Minnesota, June 1997.

[3] Deligne, P.: *Valeurs de fonctions L et périodes d'intégrales.* Proc. of Symp. in Pure Math. **33**(1979), part 2, pp. 313-346.

[4] Deligne, P., Milne, J. S., Ogus, A., Shih, K.-Y.: *Hodge Cycles, Motives, and Shimura Varieties*, Lecture Notes in Math. **900**, Springer Verlag, New York 1982.

[5] Farrell, F. T.: *An extension of Tate cohomology to a class of infinite groups.* J. Pure Appl. Algebra **10**(1977)153-161.

[6] Galovich, S. and Rosen, M.: *Units and class groups in cyclotomic function fields.* J. Number Theory **14**(1982)156-184

[7] Goss, D.: *Basic Structures of Function Field Arithmetic.* Ergebnisse der Mathematik und ihrer Grenzgebiete, 3. Folge, Vol. 35, Springer Berlin Heidelberg 1996.

[8] Hayes, D. R.: *Analytic class number formulas in function fields.* Invent. Math. **65**(1981)49-69

[9] Hayes, D.: *A brief introduction to Drinfeld modules.* in: The arithmetic of function fields (Goss, D., Hayes, D. R., Rosen, M. I., eds.), pp. 1-32, W. de Gruyter, Berlin New York 1992

[10] Kubert, D. S.: *The universal ordinary distribution.* Bull. Soc. Math. France **107**(1979), 179-202.

[11] Kubert, D. S.: *The $\mathbf{Z}/2\mathbf{Z}$-cohomology of the universal ordinary distribution.* Bull. Soc. Math. France **107**(1979)203-224

[12] Lang, S.: *Cyclotomic Fields I and II.* Graduate Texts in Mathematics **121** Springer-Verlag, New York 1990.

[13] Mislin, G.: *Tate cohomology for arbitrary groups via satellites* Topology and its Applications **56**(1994)293-300

[14] Sinha, S. K.: *Deligne's reciprocity for function fields.* J. Number Theory **63**(1997)65-88

[15] Sinnott, W.: *On the Stickelberger ideal and the circular units of a cylotomic field.* Ann. Math. **108**(1978)107-134

[16] Yin, L. S.: *Index-class number formulas over global function fields.* Compositio Math. (to appear)